\documentclass[11pt]{article}
\usepackage{amscd}
\usepackage{amssymb}
\usepackage{amsmath}
\usepackage{amsthm}
\usepackage{latexsym}
\usepackage[dvips]{graphics}
\begin{document}
\title{On the geometric simple connectivity of open manifolds\\
{\small \em IMRN 2004, to appear}}
\author{
\begin{tabular}{cc}
 Louis Funar\thanks{Partially supported by a Canon grant}
 & Siddhartha
Gadgil\thanks{{\em Current address:} 
Stat-Math Unit, Indian Statistical Institute,
8th Mile, Bangalore 560059, India}
\thanks{Partially supported by a Sloan Dissertation Fellowship}\\
\small \em Institut Fourier BP 74 &\small \em Department of
Mathematics, \\ \small \em University of Grenoble I &\small \em SUNY
at Stony Brook\\ \small \em 38402 Saint-Martin-d'H\`eres cedex, France
&\small \em Stony Brook, NY 11794-3651, USA\\ \small \em e-mail: {\tt
funar@fourier.ujf-grenoble.fr} & \small \em e-mail: {\tt
gadgil@math.sunysb.edu } \\
\end{tabular}
}

\maketitle

\newcommand{\pinf}{\pi_1^{\infty}}

\theoremstyle{definition}
\newtheorem{definition}{Definition}[section]

\theoremstyle{plain}
\newtheorem{theorem}{Theorem}[section]
\newtheorem{lemma}[theorem]{Lemma}
\newtheorem{proposition}[theorem]{Proposition}
\newtheorem{conjecture}[theorem]{Conjecture}
\newtheorem{corollary}[theorem]{Corollary}

\theoremstyle{remark}
\newtheorem{remark}[theorem]{Remark}

\def \onto{\twoheadrightarrow}
\def \del{\partial}
\def \incl{\hookrightarrow}
\def \Incl{\hookleftarrow}
\def \inte{\overset{\circ}}
\def \P{\mathfrak P}
\def \H{{\left[\begin{smallmatrix} 0 & 1\\ 1 & 0\end{smallmatrix}\right]}}
\def \Z{mathbb Z}
\def \sm{\setminus}


\begin{abstract}
A manifold is said to be
geometrically simply connected if it has a proper handle decomposition
without 1-handles. By the work of Smale, for compact manifolds of
dimension at least five, this is equivalent to simple-connectivity.
We prove that there exists an obstruction to an open simply connected
 $n$-manifold of dimension $n\geq 5$ being geometrically simply
 connected. In particular, for each $n\geq 4$ there exist uncountably many 
simply connected $n$-manifolds which are not  geometrically 
simply connected.  
We also prove that for $n\neq 4$ an $n$-manifold proper homotopy
equivalent to a weakly geometrically simply connected 
polyhedron is geometrically simply connected (for $n=4$ it is 
only end compressible). We analyze further the case $n=4$ and 
Po\'enaru's conjecture. 
\end{abstract}


end compressibility, Casson finiteness.
 

\vspace{0.5cm}
\noindent 
Since the proof of the $h$-cobordism and $s$-cobordism theorems for
manifolds whose dimension is at least $5$, a central method in
topology has been to simplify handle-decompositions as much as
possible. In the case of compact manifolds, with dimension at least
$5$, the only obstructions to canceling handles turn out to be
well-understood algebraic obstructions, namely the fundamental group,
homology groups and the Whitehead torsion.

\vspace{0.2cm}
\noindent 
In particular, a compact manifold $M$ whose dimension is at least $5$
has a handle-decomposition without $1$-handles iff it is
simply-connected. Here we study when open manifolds have a
handle-decomposition without $1$-handles.

\vspace{0.2cm}
\noindent 
In the case of $3$-manifolds, there is a classical obstruction to a
simply-connected open manifold having a handle-decomposition without
$1$-handles, namely the manifold has to be \emph{simply connected at
infinity}. However in higher dimensions there are several open
manifolds having a handle-decomposition without $1$-handles that are
not simply-connected at infinity. Thus, there is no classical
obstruction, besides the fundamental group (which is far too weak) 
and the $\pi_1$-stability at infinity (which is stronger), to an open manifold with
dimension at least $5$ having a handle-decomposition without
$1$-handles. We show here that there are (algebraic) obstructions,
which can moreover be seen to be non-trivial in concrete examples.

\vspace{0.2cm}
\noindent  
We show that for open manifolds, there is an algebraic condition,
which we call \emph{end-compressibility}, equivalent to the existence
of a proper handle-decomposition without $1$-handles. This is a
condition on the behaviour at infinity of the
manifold. End-compressibility is defined in terms of fundamental
groups related to an exhaustion of the open manifold by compact
submanifolds. In some sense it is the counterpart of Siebenmann's
obstruction to finding a boundary of an open manifold, although our
requirements are considerably weaker.

\vspace{0.2cm}
\noindent 
End-compressibility in turn implies a series of conditions, analogous
to the lower central series (the first of which is always
satisfied). Using these, we show that if $W$ is a Whitehead-type
manifold and $M$ a compact (simply-connected) manifold, $W\times M$
does not have a handle-decomposition without $1$-handles.

\vspace{0.2cm}
\noindent 
In the case of manifolds of dimension $3$ and $4$, even the compact
case is subtle. For $3$-manifolds, the corresponding statement in the
compact case is equivalent to the Poincar\'e conjecture, and
irreducible open simply-connected $3$-manifolds have a
handle-decomposition without 1-handles iff they are simply-connected
at infinity.  For a compact, contractible $4$-manifold $M$, (an
extension of) an argument of Casson shows that if $\pi_1(\del M)$ has
a finite quotient, then any handle-decomposition of $M$ has
$1$-handles. 

\vspace{0.2cm}
\noindent 
We study possible handle-decompositions of the interior of $M$ and
find a standard form for this. As a consequence, if a certain
finiteness conjecture, generalizing a classical conjecture about
links, is true, then the interior of $M$ also has no
handle-decomposition without $1$-handles.

\vspace{0.2cm}
\noindent 
More explicitly, if all the Massey products of a link $L$ in $S^3$
vanish, then it is conjectured that the link $L$ is a sub-link of a
homology boundary link. In our case, we have a link $L$ in a manifold
$K$, and a degree-one map from $K$ to $\del M$, so that the image of
$L$ in $\pi_1(\del M)$ consists of homotopically trivial curves. The
notions of {\em vanishing Massey products} and {\em homology boundary
links} generalize to versions relative to $\pi_1(M)$ in this
situation. The conjecture mentioned in the above paragraph is that
{\em links with vanishing Massey products \emph{(relative to
$\pi_1(M)$)} are sub-links of homology boundary links} (relative to
$\pi_1(M)$).

\tableofcontents

\section{Statements of the results}

\noindent The problem we address in this paper is whether 1-handles are necessary
in a handle decomposition of a simply connected manifold. 
Moreover we investigate when it is possible to kill 1-handles 
within the proper homotopy type of a given open manifold.
 
\vspace{0.2cm}
\noindent The relation between algebraic connectivity and geometric connectivity
(in various forms) was explored first by E.C.~Zeeman (see \cite{Ze}) in
connection with the Poincar\'e conjecture.  Zeeman's definition of the
geometric $k$-connectivity of a manifold amounts to asking that any
$k$-dimensional compact can be engulfed in a ball.  His main result
was the equivalence of algebraic $k$-connectivity and geometric
$k$-connectivity for $n$-manifolds, under the condition $k\leq
n-3$. Notice that it makes no difference whether one considers open and
compact manifolds. 

\vspace{0.2cm}
\noindent Later C.T.C.~Wall (\cite{Wall}) introduced another concept of geometric
connectivity using handle theory which was further developed by
V.~Po\'enaru in his work around the Poincar\'e conjecture.  A similar
equivalence between the geometric and algebraic connectivities holds
\emph{in the compact case} but this time one has to replace the
previous codimension condition by $k\leq n-4$. In this respect all
results in low codimension are hard results. There is also a
non-compact version of this definition which we can state precisely as
follows:

\begin{definition}
A (possibly non-compact) manifold, which might have nonempty boundary, is 
{\em geometrically $k$-connected} (abbreviated. {\em g.$k$.c.}) if there 
exists a proper handlebody decomposition without $j$-handles, for 
$1\leq j\leq k$. 
\end{definition}

\noindent One should emphasize that now the compact and non-compact situations
are no longer the same. The geometric connectivity is a consequence of
the algebraic connectivity only under additional hypotheses concerning
the ends. The purpose of this paper is to partially characterize these
additional conditions.

\begin{remark}
 Handle decomposition are known to exist for all manifolds in the 
topological, PL and smooth settings, except in the case of topological 
4-manifolds. In the latter case the existence of a handlebody 
decomposition is equivalent to that of a PL (or smooth) structure. 
However in the open case such a smooth structure always exist (in dimension 4).
Although most results below can be restated and proved for other 
categories, we will restrict ourselves to considering PL manifolds 
and handle decompositions in the sequel. 
\end{remark}
 
\noindent We will be mainly concerned with geometric simple connectivity (abbreviated. 
g.s.c.) in the sequel. 
A related concept, relevant only in the non-compact case, is:

\begin{definition}
A (possibly non-compact) polyhedron 
$P$ is {\em weakly geometrically simply connected} 
(abbreviated {\em w.g.s.c.}) if $P=\cup_{j=1}^{\infty} K_j$, where 
$K_1\subset K_2\subset ...\subset K_j\subset...$ 
is an exhaustion by compact connected sub-polyhedra with $\pi_1(K_j)=0$.
Alternatively, any compact subspace is contained in a simply connected 
sub-polyhedron. 
\end{definition}
\noindent Notice that a w.g.s.c. polyhedron is simply connected. 

\begin{remark}
\begin{itemize}
\item The w.g.s.c. spaces with which we will be concerned in the sequel are
usually manifolds.  Similar definitions can be given in the case of
topological (respectively smooth) manifolds where we require the
exhaustions to be by topological (respectively smooth) submanifolds.
All results below hold true for this setting too (provided
handlebodies exists) except those concerning Dehn exhaustibility,
since the later is essentially a PL concept.
\item The w.g.s.c. is much more flexible than the g.s.c., the latter
making sense only for manifolds, and enables us to work within the
realm of polyhedra. However one can easily show (see below) that
w.g.s.c.  and g.s.c. are equivalent for non-compact manifolds of
dimension different from $4$ (under the additional irreducibility
assumption for dimension $3$). Its invariance under proper homotopy
equivalences expresses the persistence of a geometric property (not
being g.s.c.) with respect to some higher dimensional manipulations
(as taking the product with a ball) of open manifolds.
\end{itemize}

\end{remark}

\noindent The first result of this paper is (see Theorem \ref{wgsc} and 
Proposition \ref{comprestrivial}):

\begin{theorem}\label{T:main}
If an open $n$-manifold is w.g.s.c., then it is end
compressible. Conversely, in dimension $n\neq 4$, if an open, simply
connected  manifold
is end compressible, then it is g.s.c. (one has to assume irreducibility
for $n=3$).
\end{theorem}

\begin{remark}
A similar result holds more generally for non-compact manifolds 
with boundary, with the appropriate definition of end compressibility. 
\end{remark}

\noindent 
End compressibility (see Definition \ref{endcompres}) is an algebraic
condition which is defined in terms of the fundamental groups of the
submanifolds which form an exhaustion.  Notice that end compressibility
is weaker than simple connectivity at infinity for $n > 3$.

\begin{remark}
The result above should be compared with Siebenmann's obstructions to
finding a boundary for an open manifold of dimension greater than 5
(see \cite{Si} and \cite{HR} for a thorough discussion of this and
related topics). The intermediary result permitting to kill 1-handles 
in this framework is theorem 3.10 p.16 from \cite{Si}: Let $W^n$ be an
open smooth $n$-manifold with  $n\geq 5$ and ${\cal E}$  an isolated
end. Assume that the end ${\cal E}$ has stable  $\pi_1$ and its
$\pi_1({\cal E})$ is finitely presented.  Then there exists
arbitrarily small 1-neighborhoods of ${\cal E}$ i.e. connected
submanifolds $V^n\subset W^n$ whose complements have compact closure, 
having compact connected boundary 
$\partial V^n$ such that $\pi_1({\cal E})\to \pi_1(V)$ and 
$\pi_1(\partial V)\to \pi_1(V)$ are isomorphisms. 
It is easy to see that this implies that the 1-neighborhood $V^n$ 
is g.s.c. This is the principal step towards canceling 
the handles of $W^n$ hence obtaining a collar. 
One notices that the hypotheses in Siebenmann's theorem are stronger 
than the end compressibility but the conclusion is stronger too.
In particular  an arbitrary w.g.s.c. manifold need not have a well-defined
fundamental group at infinity, as is the case for $\pi_1$-stable
ends. However we think that the relationship between the $\pi_1$-(semi)stability of ends 
and end compressibility would deserve further investigation. 

\vspace{0.2cm}
\noindent 
The full power of the $\pi_1$-stability is used to cancel 
more than 1-handles.  Actually L.~Siebenmann considered tame ends, which means that ${\mathcal
  E}$ is $\pi_1$-stable and it has arbitrarily small neighborhoods
which are finitely dominated. The tameness condition is strong enough
to insure (see theorem 4.5 of \cite{Si}) that all $k$-handles can be
canceled for $k\leq n-3$. One more  obstruction (the end obstruction) is actually needed
in order to be able to cancel the $(n-2)$-handles (which it turns out to imply 
the existence of a collar). There exist tame ends which are not collared (i.e. with non-vanishing
end obstruction), as well as $\pi_1$-stable ends with finitely
presented $\pi_1({\mathcal E})$ which are not tame. 
Thus  the obstructions for killing properly the  handles of index  $1\leq \lambda
\leq k$ should be weaker than the tameness of the end for $k\leq n-3$
and must coincide with Siebenmann's for $k=n-2$. 
\end{remark}

\begin{remark}
It might be worthy to compare our approach with the results from 
\cite{HR}. First the w.g.s.c. is the analogue of the 
reverse collaring. According
to (\cite{HR}, Prop.8.5, p.93) a space $W$ is reverse collared if it
has an exhaustion by compacts $K_j$ for which the inclusions 
$K_j\to W$ are homotopy equivalences (while the w.g.s.c. asks only
that these inclusions be 1-connected), and hence a simply connected 
reverse collared space is w.g.s.c. 

\vspace{0.2cm}
\noindent 
The right extension of the w.g.s.c. to non-simply connected spaces, 
which is suitable for applications
to 3-manifolds, is the Tucker property (see \cite{BT,M}), which is
also a proper homotopy invariant and can be formulated as a group
theoretical property for coverings. This is again weaker than the 
reverse tameness/collaring (see \cite{HR}, Prop.8.9 and Prop.11.13).
Moreover there is  a big difference between the
extended theory and the present one: 
while the w.g.s.c. is the property of having no extra 1-handles, 
the Tucker property expresses
the fact that some handlebody decomposition needs {\em only finitely
many} 1-handles, without any control on their number. 

\vspace{0.2cm}
\noindent 
In this respect our first result is a sharpening of the theory of
reverse collared manifolds specific to the realm of simply connected
spaces. 

\vspace{0.2cm}
\noindent 
The g.s.c. is mostly  interesting  in low dimensions, for
instance in dimension 3 where it implies the simple connectivity at
infinity. However its importance relies on its proper 
homotopy invariance, which has been
discovered in a particular form by V.~Po\'enaru (see also \cite{Fun2,FT}), 
enabling us to transform the low dimensional problem
"is $W^3$ g.s.c.?" into a high dimensional one, for example
"is $W^3\times D^n$ w.g.s.c.?". We provide in this paper a 
criterion permitting to check the answer to the high dimensional
question, in terms of an arbitrary exhaustion by compact
submanifolds. This criterion is expressed algebraically as the end
compressibility of the manifold and is closer to the forward 
tameness from \cite{HR}, rather than the reverse tameness.   
\end{remark}

\begin{remark}
If $W^k$ is compact and simply-connected then the product $W^k\times
D^n$ with a closed $n$-disk is g.s.c. if
$n+k\geq 5$.  However there exist non-compact $n$-manifolds with
boundary which are simply connected but not end compressible (hence not
w.g.s.c.) in any dimension $n$, for instance $W^3\times M^n$ where
$\pinf W^3\neq 0$. Notice that $W^k\times {\rm int}(D^n)$ is g.s.c. 
for $n\geq 1$ since $\pinf(W^k\times {\rm int}(D^n))=0$. 
\end{remark}

\noindent 
We will also prove (see Theorem \ref{uncountable}):
\begin{theorem}
There exist uncountably many open contractible $n$-manifolds for 
any $n\geq 4$ which are not w.g.s.c.  
\end{theorem}

\noindent 
The original motivation for this paper was to try to kill 1-handles of
open 3-manifolds at least {\em stably} (i.e. after stabilizing the
$3$-manifold). The meaning of the word stably in \cite{Po1}, where
such results first arose, is to do so at the expense of taking products with
some high dimensional compact ball.  This was extended in
\cite{Fun2,FT} by allowing the 3-manifold to be replaced by any other
polyhedron having the same proper homotopy type.  The analogous result is
true for $n\geq 5$ (for $n=4$ only a weaker
statement holds true (see Theorem \ref{proper}):

\begin{theorem}
If a non-compact manifold of dimension  $n\neq 4$ 
is proper homotopically dominated by  
a w.g.s.c. polyhedron, then it is w.g.s.c.
A non-compact $4$-manifold  proper homotopically 
dominated by  a w.g.s.c. polyhedron is end compressible.
\end{theorem}

\noindent Since proper homotopy equivalence implies proper homotopy
domination we obtain: 
\begin{corollary}
If a non-compact manifold of dimension $n\neq 4$ 
is proper homotopy equivalent to a w.g.s.c.
polyhedron, then it is w.g.s.c.
\end{corollary}

\begin{remark}
This criterion is an essential ingredient in Po\'enaru's proof (see
\cite{PoU}) of the covering space conjecture: If $M^3$ is a closed,
irreducible, aspherical 3-manifold, then the universal covering space
of $M^3$ is ${\bf R}^3$. Further developments suggest a similar result
in higher dimensions, by replacing the simple connectivity at infinity
conclusion with the weaker w.g.s.c. We will state below a
group-theoretical conjecture abstracting this purely 3-dimensional
result.
\end{remark}

\noindent 
It is very probable that  there exist examples of open 4-manifolds 
which are not w.g.s.c., but their products with a closed ball 
are w.g.s.c. Thus in some sense the previous result is sharp.  

\vspace{0.2cm}
\noindent 
The dimension 4 deserves special attention also because one expects that
the w.g.s.c. and the g.s.c.  are not equivalent. 
Specifically V.~Po\'enaru conjectured that:

\begin{conjecture}[Po\'enaru Conjecture]
If the interior of a compact contractible 4-manifold with boundary 
a homology sphere is g.s.c. then the compact 4-manifold is also g.s.c. 
\end{conjecture}
\begin{remark}
\begin{itemize}
\item A consequence of this conjecture, for the  particular case
of the product of a homotopy 3-disk $\Delta^3$ with an interval, is  
the Poincar\'e conjecture in dimension 3. This  follows from the
two results announced by V.~Po\'enaru:  

{\bf Theorem}: If $\Sigma^3$ is a homotopy 3-sphere such that 
$\Sigma^3\times [0,1]$ is g.s.c. then $\Sigma^3$ is g.s.c. (hence
standard). 

{\bf Theorem}: If $\Delta^3$ is a homotopy 3-disk then 
${\rm int}(\Delta^3\times [0,1]\sharp_{\infty} S^2\times D^2)$ is g.s.c.,
where $\sharp$ denotes the boundary connected sum. 
\item The differentiable Poincar\'e conjecture in dimension 4 is
  widely believed to be false. One reasonable reformulation of it 
would be the following: {\em A smooth homotopy 4-sphere (equivalently, 
homeomorphic to $S^4$) that is g.s.c. should be diffeomorphic to 
$S^4$}. 
\item The two conjectures above (Po\'enaru's and the reformulated Poincar\'e)
  conjectures imply also the smooth Schoenflies conjecture in
  dimension 4, which states that a $3$-sphere smoothly embedded in $S^4$
bounds a smoothly embedded $4$-ball. In fact by a celebrated result of 
B.~Mazur any such Schoenflies ball has its interior diffeomorphic to
${\bf R}^4$, hence g.s.c.
      
\end{itemize}
\end{remark}

\noindent 
An immediate corollary would be that the interior of a 
Po\'enaru-Mazur 4-manifold may be w.g.s.c. but not g.s.c., because 
some (compact) Po\'enaru-Mazur 4-manifolds are known to be not g.s.c.
(the geometrization conjecture implies this statement for all 
4-manifolds whose boundary 
is not a homotopy sphere). The proof is due to A.~Casson 
and it was based on partial positive solutions to the
following algebraic conjecture \cite[p.117]{Ke} and \cite[p.403]{MKS}.

\begin{conjecture}[Kervaire Conjecture]\label{kerv} Suppose one adds an equal
number of generators $\alpha_1,\dots,\alpha_n$ and relations
$r_1,\dots,r_n$ to a non-trivial group $G$, then the group
$\frac{G*\langle\alpha_1,\dots,\alpha_n\rangle}
{\langle\langle r_1,\dots,r_n\rangle\rangle}$ that one
obtains is also non-trivial.
\end{conjecture}

\noindent 
A.~Casson showed that certain $4$-manifolds $(W^4,\del W^4)$ have no handle
decompositions without $1$-handles by showing that if they did, then
$\pi_1(\del W^4)$ violates the Kervaire conjecture. 
Our aim would be to show that most contractible 4-manifolds 
are not g.s.c., and the method of the proof is to reduce 
this statement to the compact case. 
However our methods permit us to obtain only a weaker result, in which 
one shows that the interior of such a manifold cannot have 
handlebody decompositions without 1-handles, if the decomposition 
has also some additional properties (see 
Theorems \ref{embed} and \ref{imers} for precise statements): 

\begin{theorem}
Assume that we have a proper handlebody decomposition without
1-handles for the interior of a Po\'enaru-Mazur 4-manifold.  If there
exists a far away intermediary level 3-manifold $M^3$ whose homology
is represented by disjoint embedded surfaces and whose fundamental
group projects to the trivial group on the boundary, then the compact
4-manifold is also g.s.c. There always exists a collection of immersed
surfaces, which might have non-trivial intersections and
self-intersections along homologically trivial curves, that fulfills
the previous requirements.
\end{theorem}

\begin{remark}
Almost all of this paper deals with geometric 1-connectivity.  However
the results can be reformulated for higher geometric connectivities
within the same range of codimensions.
\end{remark}


\noindent 
We wish to emphasize that there is a strong group theoretical flavour
in the w.g.s.c. for universal covering spaces. In this
respect the universal covering conjecture in dimension three (see
\cite{PoU}) would be a first step in a more general program. Let us define
a finitely presented, infinite group $\Gamma$ to be {\em w.g.s.c.} if
there exists a compact polyhedron with fundamental
group $\Gamma$ whose universal covering space is
w.g.s.c. It is not hard to show that this
definition does not depend on the particular polyhedron one chooses
but only on the group.  This is part of a more general philosophy, due
to M.~Gromov, in which infinite groups are considered as geometric
objects.  This agrees with the idea that killing 1-handles of
manifolds is a group theoretical problem in topological disguise. The
authors think that the following might well be true:

\begin{conjecture}
Fundamental groups of closed aspherical manifolds 
are w.g.s.c. 
\end{conjecture}

\noindent This will be a far
reaching generalization of the three dimensional result 
announced by  V.~Po\'enaru in \cite{PoU}. 
It is worthy to note that all reasonable examples of groups 
(e.g. word hyperbolic, semi-hyperbolic, $CAT(0)$, group
extensions, one relator groups) are w.g.s.c.  
It would be interesting to find an example of a finitely presented group which fails
to be w.g.s.c. Notice that the well-known examples of M.~Davis of 
Coxeter groups which are fundamental groups 
of aspherical manifolds whose universal covering spaces are not simply
connected at infinity are actually $CAT(0)$ hence w.g.s.c. 
However one might expect a direct connection between the
semi-stability of finitely presented groups, the quasi-simple
filtrated groups (see \cite{BM}) and the
w.g.s.c. We will address these questions in a future paper.

\subsection*{Outline of the paper}

\noindent 
In section 2, we compare w.g.s.c., g.s.c and  the simple connectivity
at infinity (s.c.i.), showing that
w.g.s.c and g.s.c. are equivalent in high dimensions and presenting
some motivating examples. 

\vspace{0.2cm}
\noindent 
Section 3 contains the core of the paper,
where we introduce the algebraic conditions and prove their relation
to w.g.s.c. We then construct uncountably many Whitehead-type
manifolds in section~4, and show that there are uncountably many
manifolds that are not geometrically simply-connected.

\vspace{0.2cm}
\noindent 
In section 5, we show that end-compressibility is a proper-homotopy
invariant. Finally, in sections~6 and~7, we turn to the
$4$-dimensional case.

\vspace{0.2cm}

\noindent{\bf Acknowledgements}. We are indebted to 
David Gabai, Ross Geoghegan,
Valentin Po\'enaru, Vlad Sergiescu and  Larry Siebenmann for valuable 
discussions and suggestions. We are grateful to the referee for his
comments which permitted us to improve the clarity and 
the readability of this paper. 
Part of this work has been done when the first author
visited Tokyo Institute of Technology, which he wishes to thank for
their 
support and hospitality, and especially to Teruaki Kitano, Tomoyoshi
Yoshida and Akio Kawauchi. Part of this work was done while the second
author was supported by a Sloan Dissertation fellowship.

\section{On the g.s.c.  }
\subsection{Killing 1-handles of 3-manifolds after stabilization}

\noindent We start with some motivating remarks about the compact 3-dimensional
situation, for the sake of comparison.
\begin{definition}
The {\em geometric 1-defect} $\epsilon(M^n)$ of the compact manifold $M^n$
is $\epsilon(M^n)=\mu_1(M^n)-{\rm rank}\; \pi_1(M)$, where $\mu_1(M^n)$ is
the minimal number of 1-handles in a handlebody decomposition and
${\rm rank}\; \pi_1(M)$ is the minimal number of generators of $\pi_1(M)$.
\end{definition}

\begin{remark}
The defect (i.e., the geometric 1-defect) is always non-negative.
There exist examples (see \cite{BZ}) of 3-manifolds $M^3$ with ${\rm rank}\;
\pi_1(M^3)=2$ and defect $\epsilon(M^3)=1$. No examples of 3-manifolds
with larger defect, nor of closed 4-manifolds with
positive defect are presently known. However it's probably true that
$\epsilon(M^3\times [0,1])=0$, for all closed 3-manifolds though 
this might be difficult to settle even for the explicit examples of
M.~Boileau and H.~Zieschang.
\end{remark}

\noindent The defect is meaningless in high dimensions because of:

\begin{proposition}
For a compact manifold  $\epsilon(M^n)=0$ holds, if $n\geq 5$. 
\end{proposition}

\begin{proof}
The proof is standard. Consider a 2-complex $K^2$ associated to a
presentation of $\pi_1(M^n)$ with the minimal number $r$ of
generators. By general position there exists an embedding
$K^2\hookrightarrow M^n$ inducing an isomorphism of fundamental
groups. Then a regular neighborhood of $K^2$ in $M^n$ has a handlebody
decomposition with $r$ 1-handles. Since the complement is 1-connected
then by (\cite{Wall,Quinn}) it is g.s.c. for $n\geq 5$ and this yields
the claim.
\end{proof}
\begin{corollary}
For a closed 3-manifold $M^3$ one has $\epsilon(M^3\times D^2)=0$. 
\end{corollary}
\begin{remark}
As a consequence if $\Sigma^3$ is a homotopy 3-sphere then
$\Sigma^3\times D^2$ is g.s.c. Results of Mazur (\cite{Mazur},
improved by Milnor in dimension 3) show that $\Sigma^3\times
D^3=S^3\times D^3$, but it is still unknown whether $\Sigma^3\times
D^2=S^3\times D^2$ holds. An earlier result of Po\'enaru states that
$(\Sigma^3-n D^3)\times D^2=(S^3-n D^3)\times D^2$ for some $n\geq 1$.
More recently, Po\'enaru's program reduced the Poincar\'e Conjecture
to the g.s.c. of $\Sigma^3\times [0,1]$.
\end{remark}

\subsection{$\pinf$ and g.s.c.}

\noindent The remarks which follow are intended to (partially) clarify the
relation between g.s.c. and the simple connectivity at infinity 
(which will be abbreviated as {\em s.c.i.} in the sequel), in general. 

\vspace{0.2cm}
\noindent Recall that a space $X$ is {\em s.c.i.} (and one also writes $\pinf(X)=0$)
if for any compact $K\subset X$ there exists a larger compact 
$K\subset L\subset X$ having the property that, any loop 
$l\subset X-L$ is null homotopic within $X-K$.  This is an important
tameness condition for the ends of the space.   
The following result was proved in (\cite{PoTa}, Thm. 1): 
\begin{proposition}
Let $W^n$ be an open simply connected $n$-manifold of dimension 
$n\geq 5$. If $\pinf(W^n)=0$ then  
$W^n$ is g.s.c.
\end{proposition}

\begin{remark}
The converse fails as the following examples show. 
Namely, for any $n\geq 5$ there exist open $n$-manifolds 
$W^n$ which are geometrically $(n-4)$-connected but  
$\pinf(W^n)\neq 0$. 

\vspace{0.2cm}
\noindent 
There exist compact contractible $n$-manifolds $M^n$ with
$\pi_1(\partial M^n)\neq 0$, for any $n\geq 4$ (see
\cite{Mazur0,Po,Glaser}).  Since $k$-connected compact $n$-manifolds
are geometrically $k$-connected if $k\leq n-4$ (see
\cite{Quinn,Wall}), these manifolds are geometrically
$(n-4)$-connected.  Let us consider now $W^n={\rm int}(M^n)$, which is
diffeomorphic to $M^n\cup_{\partial M^n\cong \partial M^n\times\{0\}}
\partial M^n \times[0,1)$.  Any Morse function on $M^n$ extends over
${\rm int}(M^n)$ to a proper one which has no critical points in the open
collar $\partial M^n \times[0,1)$, hence ${\rm int}(M^n)$ is also
geometrically $(n-4)$-connected.  On the other hand $\pi_1(\partial
M^n)\neq 0$ implies $\pinf(W^n)\neq 0$. 

\end{remark}

\noindent However the following partial converse holds:

\begin{proposition}\label{cute}
Let $W^n$ be a non-compact simply connected $n$-manifold which 
has a proper handlebody decomposition
\begin{enumerate}
\item without $1-$ or $(n-2)-$handles, or
\item without $(n-1)-$ or $(n-2)-$handles.
\end{enumerate}
Then $\pinf(W^n)=0$.
\end{proposition}

\begin{remark}
When $n=3$ this simply says that 1-handles are necessary unless 
$\pinf(W^3)=0$.
\end{remark}


\noindent {\em Proof of Proposition \ref{cute}}.
Consider the handlebody decomposition $W^n = B^n \cup_{j=1}^{\infty}
h_j^{{i_j}}$, where $h_j^{{i_j}}$ is an ${i_j}$-handle ($B^n=h_0^0$).
Set $X_m = B^n\cup_{j=1}^{m}h_j^{{i_j}}$, for $m\geq 0$.  Assume that
this decomposition has no $1-$ nor $(n-2)-$handles.  Since there are
no 1-handles it follows that $\pi_1(X_j)=0$ for any $j$ (it is only
here one uses the g.s.c.).

\begin{lemma}
If $X^n$ is a compact simply connected  $n$-manifold having a handlebody 
decomposition without $(n-2)-$handles then  $\pi_1(\partial X^n)=0$.
\end{lemma}
\begin{proof}
Reversing the handlebody decomposition of $X^n$ one finds a decomposition 
from $\partial X^n$ without 2-handles. One slides the handles to be attached 
in increasing order of their indices.
Using Van Kampen Theorem it follows that 
$\pi_1(X^n)=\pi_1(\partial X^n)* {\bf F}(r)$, where $r$ is the number of 
1-handles, and thus $\pi_1(\partial X^n)=0$. \end{proof}

\begin{lemma}
If the compact submanifolds $...\subset X_m\subset X_{m+1}\subset ...$ 
exhausting the simply connected manifold $W^n$
satisfy  $\pi_1(\partial X_m)=0$, for all $m$, then $\pinf(W^n)=0$.
\end{lemma}

\begin{proof}
For $n=3$ this is clear. Thus we suppose $n\geq 4$.  For any compact
$K\subset W^n$ choose some $X_m\supset K$ such that $\partial X_m \cap
K=\emptyset$.  Consider a loop $l\subset W^n-X_m$. Then $l$ bounds an
immersed (for $n\geq 5$ embedded) 2-disk $\delta^2$ in $W^n$. We can
assume that $\delta^2$ is transversal to $\partial X_m$. Thus it
intersects $\partial X_m$ along a collection of circles
$l_1,...,l_p\subset\partial X_m$.  Since $\pi_1(\partial X_m)=0$ one is able to
cap off the loops $l_j$ by some immersed 2-disks $\delta_j\subset
\partial X_m$. Excising the subsurface $\delta^2\cap X_m$ and
replacing it by the disks $\delta_j$ one obtains an immersed 2-disk
bounding $l$ in $W^n-K$.  \end{proof} 

\noindent This proves the first case of the Proposition \ref{cute}.
In order to prove the second case choose some connected compact subset
$K\subset W^n$. By compactness there exists $k$ such that $K\subset
X_k$.  Let $r$ be large enough (this exists by the properness)
such that any handle $h_p^{i_p}$ whose attaching zone touches the
lateral surface of one of the handles $h_1^{i_1}, h_2^{i_2},...,
h_k^{i_k}$ satisfies $p\leq r$.  The following claim will prove the
Proposition \ref{cute}:

\begin{proposition}\label{gicu}
Any loop $l$ in $W^n-X_r$ is null-homotopic in $W^n-K$.
\end{proposition}

\begin{proof} Actually the following more general engulfing  result holds:

\begin{proposition}\label{engulf}
If $C^2$ is a 2-dimensional polyhedron whose boundary 
$\partial C^2$ is contained in $W^n-X_r$ then there exists 
an isotopy of $W^n$ (with compact support), fixing $\partial C^2$ 
and moving $C^2$ into  $W^n-K$. 
\end{proposition}

\begin{proof}
Suppose that  $C^2 \subset X_m$. 
One reverses the handlebody decomposition of $X_m$ and obtains 
a decomposition from $\partial X_m$ without 1- or 2-handles. 
 Assume that we can move $C^2$ such that it misses the last 
$j\leq r-1$ handles. By general position there exists an isotopy 
(fixing the last $j$ handles) making 
$C^2$ disjoint of the co-core ball of the $(j-1)$-th handle, since 
the co-core disk has dimension at most $n-3$.
The uniqueness of the regular neighborhood implies that 
we can move $C^2$ out of the $(j-1)$-th handle (see e.g. \cite{RS}), by 
an isotopy which is identity on the last $j$ handles. 
This proves  the Proposition \ref{engulf}. 
\end{proof}
\noindent This yields the result of Proposition \ref{gicu} 
by taking for $C^2$ any 2-disk
parameterizing a null homotopy of $l$.
\end{proof}

\subsection{G.s.c. and w.g.s.c.}
\begin{proposition}
The non-compact manifold $W^n$ ($n\neq 4$), which one supposes to be
irreducible if $n=3$, is w.g.s.c. if and only if it is g.s.c.
\end{proposition}

\begin{proof}
For $n=3$ it is well-known that g.s.c. is equivalent to w.g.s.c. 
which is also equivalent to s.c.i. if the manifold is 
irreducible. For $n\geq 5$ 
this is a consequence of Wall's result stating the 
equivalence of g.s.c. and simple connectivity 
in the compact case (see \cite{Wall}).
If  $W^n$ is w.g.s.c. then it has an exhaustion by compact 
simply connected sub-manifolds $M_j$ (by taking suitable regular 
neighborhoods of the polyhedra). One can also refine the exhaustion 
such that the boundaries are disjoint. Then the pairs 
$({\rm cl}(M_{j+1}-M_j), \partial M_j)$ are 1-connected, hence  (\cite{Wall}) 
they have a handlebody decomposition without 1-handles. 
Gluing together these intermediary decompositions we obtain a 
proper handlebody decomposition as claimed.
\end{proof}

\section{W.g.s.c. and end compressibility}

\noindent 
In this section, we show that w.g.s.c. is equivalent to an algebraic
condition which we call \emph{end compressibility}. This in turn
implies infinitely many conditions, \emph{end $k$-compressibility} for
ordinals $k$, and is equivalent to all these plus a finiteness
condition.

\vspace{0.2cm}
\noindent 
Using the above, we give explicit examples of open manifolds that are
not w.g.s.c.

\subsection{Algebraic preliminaries}

\noindent 
In this section we introduce various algebraic notions of
{\em compressibility} and study the relations between these. This will be
applied in a topological context in subsequent sections, where
{\em compressibility} corresponds to being able to attach {\em enough} two
handles, and {\em stable-compressibility} refers to the same after
possibly attaching some $1$-handles.

\begin{definition}
A pair $(\varphi:A\to B, \; \psi:A\to C)$ of group morphisms  
is {\em strongly  compressible} if $\varphi(\ker \psi)=\varphi(A)$.
\end{definition}

\begin{remark}
Strong compressibility is symmetric in the arguments $(\varphi,\psi)$
i.e. $\varphi(\ker \psi)=\varphi(A)$ is equivalent to
$\psi(\ker\varphi)=\psi(A)$. The proof is an elementary diagram chase.
\end{remark}

\begin{definition}
A pair $(\varphi:A\to B, \; \psi:A\to C)$ of group morphisms,  
is {\em stably compressible} if there exists 
some free group ${\bf F}(r)$ on finitely many generators, and 
a morphism $\beta:{\bf F}(r)\to C$, such that the pair 
 $(\varphi * 1_{{\bf F}(r)}:A* {\bf F}(r)\to B*{\bf F}(r), \; 
\psi* \beta:A *{\bf F}(r)\to C)$ is strongly compressible.
\end{definition}

\noindent 
We shall see that stable-incompressibility implies infinitely many
conditions, indexed by the ordinals, on the pair of morphisms. We
first define a series of groups (analogous to the lower central
series).

\begin{definition}
Consider a fixed pair $(\varphi:A\to B, \; \psi:A\to C)$ of group
morphisms.  We define inductively a subgroup $G_{\alpha}\subset C$ for
any ordinal $\alpha$. Set $G_0=C$.  If $G_{\alpha}$ is defined for
every $\alpha <\beta$ (i.e. $\beta$ is a limit ordinal) then set
$G_{\beta}=\cap_{\alpha < \beta} G_{\alpha}$. Further set
$G_{\alpha+1}={\cal N}(\psi(\ker \varphi), G_{\alpha})\triangleleft
G_{\alpha}$ for any other ordinal, where ${\cal N}(K, G)$ is the
smallest normal group containing $K$ in $G$.  
\end{definition}

\noindent 
The groups $G_{\alpha}$ form a decreasing sequence of subgroups of
$C$. Using Zorn's lemma there exists an infimum of the lattice of
groups $G_{\alpha}$, ordered by the inclusion, which we denote by
$G_{\infty}=\cap_{\alpha}G_{\alpha}$ (over all ordinals $\alpha$).

\begin{definition}
The pair $(\varphi:A\to B,\psi:A\to C)$ is said to be
{\em $\alpha$-compressible} if $\psi(A) \subset G_{\alpha}$ (where $\alpha$
is an ordinal or $\infty$).
\end{definition}

\begin{lemma}
Given a subgroup $L\subset C$ there exists a maximal 
subgroup $\Gamma=\Gamma(L,C)$ of $C$ so that 
$L\subset {\cal N}(L,\Gamma)=\Gamma$. 
\end{lemma}

\begin{proof}
There exists at least one group $\Gamma$, for instance $\Gamma=L$.
Further if $\Gamma$ and $\Gamma'$ verify the condition
${\cal N}(L,\Gamma)=\Gamma$, then their product $\Gamma\Gamma'$ does.  Thus,
Zorn's Lemma says that a maximal element for the lattice of subgroups
verifying this property (the order relation  is the inclusion) exists.
\end{proof}

\begin{lemma}
We have $\Gamma(L, C)=G_{\infty}$, where 
$L=\psi(\ker \varphi)$. 
\end{lemma}

\begin{proof}
First, $G_{\infty}$ satisfies the condition ${\cal N}(L,\Gamma)=\Gamma$
otherwise the minimality will be contradicted.  Pick an arbitrary
$\Gamma$ satisfying this condition.  If $\Gamma\subset G_{\alpha}$ it
follows that $\Gamma={\cal N}(L,\Gamma) \subset {\cal N}(L,
G_{\alpha})=G_{\alpha+1}$, hence by a transfinite induction we derive
our claim. 
\end{proof}

\begin{definition}
One says that $K$ is {\em full} in $\Gamma$ if ${\cal N}(K,\Gamma)=\Gamma$. 
If we have a pair and $\psi(\ker \varphi)$ is full in $\Gamma$, then  
we call $\Gamma$ {\em admissible}. 
\end{definition}

\begin{remark}
If  $\Gamma$ is admissible 
then $\psi(\ker\varphi)\subset G_{\infty}$, since $G_{\infty}$ is the largest 
group with this property.
\end{remark}

\begin{proposition} \label{compres}
If the pair $(\varphi:A\to B, \; \psi:A\to C)$, 
where $A, B$ and $C$ are finitely generated and $\varphi(A)$ 
is finitely presented   
is stably compressible then it is $\infty$-compressible
and there exists a subgroup  $\Gamma\subset G_{\infty}\subset C$ which 
is normally finitely generated  within $C$  and 
such that $\psi(\ker \varphi)$ is full in $\Gamma$. 

\vspace{0.2cm}
\noindent 
Conversely, if the pair $(\varphi:A\to B, \; \psi:A\to C)$  is 
$\infty$-compressible
and there exists a finitely generated  subgroup  $\Gamma\subset G_{\infty}$  
such that $\psi(\ker \varphi)$ is full in $\Gamma$, then 
the pair is stably compressible.
\end{proposition}

\begin{proof}
We set $K=\ker\varphi$ in the sequel. We establish first:

\begin{lemma}\label{stabtoinf}
If the pair $(\varphi:A\to B, \; \psi:A\to C)$ 
is stably compressible then it is $\infty$-compressible.
\end{lemma}

\begin{proof}
We will use a transfinite recurrence with the inductive steps provided
by the next two lemmas. Set $\beta:{\bf F}(r)\to C$ for the morphism
making the pair $(\varphi *1, \psi*\beta)$ strongly compressible.
\begin{lemma}
If $\beta({\bf F}(r))\subset G_{i}$ and $\psi(A)\subset G_i$ 
then $\psi(A)\subset G_{i+1}$. 
\end{lemma}
\begin{proof}
By hypothesis $\varphi*1(\ker \psi *\beta)\supset \varphi *1 (A*{\bf
  F}(r))$. 
Alternatively, for any $b\in \varphi(A)\subset B\subset  B* {\bf
  F}(r)$ there exists 
some $x\in A * {\bf F}(r)$ such that $\varphi * 1(x)= b$ 
and $\psi * \beta(x)=1$.
One can write uniquely $x$ in normal form (see \cite{LS}, Thm.1.2., p.175) 
as $x=a_1f_1a_2f_2...a_mf_m$ where 
$a_j\in A, f_j\in {\bf F}(r)$ are non-trivial (except maybe $f_m)$.
Then 
$\varphi*1(x)=\varphi(a_1)f_1\varphi(a_2)f_2...\varphi(a_m)f_m$. 

\vspace{0.2cm}
\noindent 
Since the normal form is unique in $B* {\bf F}(r)$ one derives that
$x$ has the following property. There exists a sequence 
$p_0=1 < p_1 < ...< p_l\leq m$ of integers  for which 
\[\varphi(a_{p_j})=b_j\neq 1\in B, \mbox{ where }  
b_1b_2...b_j=b, \] 
\[\varphi(a_j)=1, \mbox{ for all } j\not\in \{p_0,p_1,...,p_l\},\] 
and 
\[ f_{p_j}f_{p_{j}+1}...f_{p_{j+1}-1}=1, \mbox{ (for all $j$, with the convention  
$p_{l+1}=m$).} \]
Furthermore $1=\psi *\beta(x)$ implies that (recall that $K=\ker\varphi$)
\[ 1\in \psi(a_1K)\beta(f_1)\psi(K)\beta(f_2)\psi(K)...\]
\[...\beta(f_{{p_1}-1})
\psi(a_{p_1}K)\beta(f_{p_1})\psi(K)\beta(f_{{p_1}+1})\psi(K)...\beta(f_m).\]
However each partial product starting at  the $p_j$-th term and ending at the 
$(p_{j+1}-1)$-th term is a product of conjugates of $\psi(K)$ by
elements from the image of $\beta$:
\[ \beta(f_{p_j})\psi(K)\beta(f_{{p_j}+1})\psi(K)...\psi(K)
\beta(f_{p_{j+1}-1})= \]
\[\prod_{i=0}^{p_{j+1}-p_j-1}
\left ( \beta\left(\prod_{k=0}^{i}f_{{p_j}+k}\right) \psi(K)
\beta\left(\prod_{k=0}^{i}f_{{p_j}+k}\right)^{-1}\right)
\subset {\cal N}(\psi(K), G_i)=G_{i+1}.\]
We used above the inclusions  $\psi(K)\subset \psi(A)\subset G_i$
and $\beta({\bf F}(r))\subset G_i$. 
Therefore 
\[ 1\in \psi(a_1K)G_{i+1}\psi(a_{p_1}K)G_{i+1} ...
\psi(a_{p_l}K)G_{i+1} =\]
\[= \psi(aK)G_{i+1}, \]
for any $a\in A$ such that $\varphi(a)=b$. This implies that 
$\psi(aK)\subset G_{i+1}$ and hence 
$\psi(A)\subset G_{i+1}$. 
\end{proof}

\begin{lemma}
If $\beta({\bf F}(r))\subset G_{i-1}$ and $\psi(A)\subset G_i$ 
then $\beta({\bf F}(r))\subset G_{i}$. 
\end{lemma}

\begin{proof}
One can use the symmetry of the algebraic compressibility 
and then the argument from the previous lemma.  
Alternatively, choose $f\in {\bf F}(r)\subset  B* {\bf
  F}(r)$ and 
some $x\in A * {\bf F}(r)$ such that $\varphi * 1(x)= f$ 
and $\psi * \beta(x)=1$. Using the normal form as above we find this 
time $1\in \beta(f){\cal N}(\psi(K), G_{i-1})$ hence 
$\beta({\bf F}(r))\subset G_i$. 
\end{proof}

\noindent Using in an alternate way the two previous lemmas one gets lemma 
\ref{stabtoinf}. 
\end{proof}

\begin{lemma}
Let $\beta:{\bf F}(r)\to C$ be a homomorphism 
such that  $(\varphi*1, \; \psi *\beta)$ 
is strongly compressible. Set $\beta({\bf F}(r))=H$. Then 
$\psi(K)$ is full in $\psi(K)H$. In particular if $A,B,C$ are 
finitely generated and $\varphi(A)$ is finitely presented then 
the subgroup $\Gamma= \psi(K)H$ is finitely generated. 
\end{lemma}

\begin{proof}
We already saw that $H\subset G_{\infty}$. 
Set $W(L; X)=\{x\mid x=\prod_i g_ix_ig_i^{-1}, \;\; 
g_i\in X, \;, x_i\in L \}$ for two subgroups $L, X\subset C$. 
The proof we used to show that $\psi(A)\subset G_{\infty}$ 
and $H\subset G_{\infty}$ actually yields 
$\psi(A)\subset W(\psi(K), H)$ and respectively 
$H\subset W(\psi(K), H)$. 
We remark now that $W(\psi(K), H)={\cal N}(\psi(K), \psi(K)H)$. 
The left inclusion is obvious. The other inclusion consists in 
writing any element $gxg^{-1}$ with $g\in \psi(K)H$, $x\in \psi(K)$ 
as a product of conjugates by elements of $H$. This might be done by
recurrence on the length of $g$, by using the following trick. 
If $g=y_1a_1y_2a_2$, $a_i\in \psi(K), y_i\in H$ then 
$gxg^{-1}=y_1 (a_1y_2a_1^{-1})(a_1a_2xa_2^{-1}a_1^{-1})
a_1y_2a_1^{-1}y_1^{-1}$. 

\vspace{0.2cm}
\noindent 
Consequently the fact that $H\subset W(\psi(K), H)$ implies 
\[W(\psi(K), H)\subset W(\psi(K), W(\psi(K), H))={\cal N}(\psi(K),
{\cal N}(\psi(K),\psi(K)H))\subset \]
\[ \subset  {\cal N}(\psi(K), \psi(K)H)=W(\psi(K), H), \]
hence all inclusions are equalities.
Also $\psi(K)H\subset {\cal N}(\psi(K), \psi(K)H)$ 
since both components $\psi(K)$ and $H$ are contained in 
${\cal N}(\psi(K), \psi(K)H)$. This shows that 
${\cal N}(\psi(K), \psi(K)H)=\psi(K)H$ 
hence  $\psi(K)$ is full in  $\psi(K)H$. 

\vspace{0.2cm}
\noindent 
We take therefore $\Gamma=\psi(K)H$. It suffices to show now that 
each of the groups $K$ and $H$ are finitely generated. 
$H$ is finitely generated since it is the image of ${\bf F}(r)$. 
Furthermore $K$ is finitely generated since $A/K=\varphi(A)$ 
is finitely presented and $A$ is finitely generated. 
The theorem of Neumann (\cite{Ba}, p.52) shows that $K$ must be 
normally finitely generated. This proves the claim.
\end{proof}

\begin{lemma}\label{admiss}
Assume that $\psi(K)$ is full in  
$\Gamma\subset G_{\infty}$, where $\Gamma$  is finitely generated. 
If the pair $(\varphi,\psi)$ is  $\infty$-compressible then 
it is stably compressible.
\end{lemma} 
\begin{proof}
Consider $r$ big enough  and a surjective homomorphism 
$\beta:{\bf F}(r)\to \Gamma$. 
This implies that $\psi *\beta (A*{\bf F}(r))=\psi(K)\Gamma=\Gamma$.  
We have to show that any $x\in \Gamma$ is in $\psi *\beta(\ker
\varphi *1)$. 

\vspace{0.2cm}
\noindent 
Recall that ${\cal N}(\psi(K), \Gamma)=\Gamma$. Set
$\;^{g}x=gxg^{-1}$. 
Then $x=\prod_i \;\;^{g_i}x_i$ can be written as a product of conjugates 
of elements $x_i\in\psi(K)$ by elements $g_i\in \Gamma$. 
Choose $f_i\in {\bf F}(r)$ so that $\beta(f_i)=g_i$ and 
$y_i\in K$ so that $\psi(y_i)=x_i$. 
Then $\psi*\beta(\prod_i f_i^{-1}y_if_i) = x$ and 
$\varphi * 1(\prod_i f_i^{-1}y_if_i) = 
 \prod_i f_i^{-1}\varphi(y_i)f_i=1$, since $y_i\in K$.
\end{proof}

\noindent Then the proposition \ref{compres} follows. 
\end{proof}

\subsection{End compressible manifolds}

\begin{definition}\label{comp}
The pair of spaces $(T',T)$ is {\em strongly compressible} (respectively 
{\em stably compressible})  
if for each component $S_j$ of $\del T$ 
one chooses a component $V_j$ of $T'-{\rm int}(T)$ such that 
$S_j\subset V_j$ so that the pair 
$(*_j\pi_1(S_j)\to \pi_1(T), \; *_j\pi_1(S_j)\to *_j\pi_1(V_j))$ is
strongly compressible (respectively stably compressible).  
The morphisms are induced by the obvious inclusions. 
Similarly, the pair  
of spaces $(T',T)$ is said to be $\alpha$-{\em compressible}  
if the pair of morphisms from above is $\alpha$-
compressible. 
Set also 
$G_{\infty}(T, T')$ for the $G_{\infty}$ group associated to 
this pair of morphisms. 
\end{definition}

\begin{remark} These morphisms are not uniquely defined and depend 
on the various choices of base points in each component. 
However the compressibility does not depend on the particular 
choice of the representative. 
\end{remark}

\begin{definition}\label{endcompres}
The open manifold $W^n$ is {\em end compressible} (respectively
{\em end $k$-compressible}) if every exhaustion of $W^n$ by compact submanifolds
$T_i^n$, such that $\pi_1(\del T_i^n)\onto\pi_1(T_i^n)$ is a surjection, has
a refinement
\[W^n=\bigcup_{i=1}^{\infty}T_i^n, \,\, T_i^n\subset {\rm int}(T_{i+1}^n), \] 
such that: 
\begin{enumerate}
\item all pairs $(T_{i+1}^n,T_i^n)$ are stably-compressible 
(respectively $k$-compressible). 
\item if $S_{i,j}^{n-1}$ denote the components of $\del T_i^n$ then the 
homomorphism  
$*_j:\pi_1(S_{i,j}^{n-1})\to \pi_1(T_i^n)$ induced by the inclusion 
is surjective. 
\item any component of $T_{i+1}^n-{\rm int}(T_i^n)$ intersect $T_i^n$ along
  precisely one component. 
\end{enumerate} 
\end{definition}

\begin{remark}
As in the case of the compressibility the condition 2 above is
independent of the homomorphism we chose, which might depend on the 
base points in each component.
\end{remark}

\begin{remark}
\begin{itemize}
\item One can ask that each connected component of 
$T_{i+1}^n-{\rm int}(T_i^n)$ has
exactly one boundary component from $\partial T_i^n$. By adding to an
arbitrary given $T_i^n$ the regular neighborhoods of arcs in
$T_{i+1}^n-{\rm int}(T_i^n)$ joining different connected component this
condition will be fulfilled.  

\item Any simply-connected manifold $W$ of dimension at least $5$ has
an exhaustion by $T_i^n$ that have the property that the natural maps
$\pi_1(\partial T_i^n)\to \pi_1(T_i^n)$ are surjective for all $i$. A
proof will be given in the next section (see lemmas \ref{T:expand} and
\ref{T:exp2}). Thus the above condition is never vacuous.
\end{itemize}
\end{remark}

\noindent 
We shall henceforth assume that exhaustions have the property that the
natural maps $\pi_1(\partial T_i^n)\to \pi_1(T_i^n)$ are surjective for
all $i$.

\begin{remark} 
The condition that the pair $(T_{i+1}^n,T_i^n)$ is stably-compressible
is implied by (and later it will be proved that it is equivalent to)
the pair of conditions
\begin{enumerate}
\item $(T_{i+1}^n,T_i^n)$ is $\infty$-compressible
\item There exists an admissible subgroup $\Gamma_i$ of
$G_{\infty}(T_i^n, T_{i+1}^n)$ which is finitely presented.
\end{enumerate}
\end{remark}

\begin{theorem}\label{wgsc}
Any  w.g.s.c. open $n$-manifold  $W^n$ is end compressible.
Conversely, for $n\neq 4$, $W^n$ is end compressible if and 
only if it is w.g.s.c.
\end{theorem}

\begin{remark}
Notice that the end compressibility of $W^3$ implies that of $W^3\times
D^2$. As a consequence of  this result for $n\geq 5$  we will derive that 
$W^3\times D^2$ is w.g.s.c. and the invariance of the w.g.s.c. under 
proper homotopies (see theorem \ref{proper}) will imply the result of 
the theorem for $n=3$.  We will restrict then for the proof to $n\geq 5$.
\end{remark}

\begin{remark}\label{indep}
It is an important issue to know whether the stable-compressibility of
{\em one particular exhaustion} implies the stable-compressibility 
of some refinement of {\em any exhaustion}. This is a
corollary of our theorem \ref{wgsc} and proposition \ref{T:wgsc}. 
In fact if $W^n$ is as above then $W^n\times D^k$ has one
stably-compressible exhaustion. Take  $n+k\geq 5$ to insure that  
$W^k\times D^n$ is w.g.s.c. and use the proposition \ref{T:wgsc}. 
In particular any product exhaustion has a stably-compressible
refinement, and the claim follows. 
\end{remark}

\subsection{Proof of Theorem \ref{wgsc}}

\noindent 
Let us consider an exhaustion $\{T_i^n\}_{i=1,\infty}$ of $W^n$ by compact 
submanifolds, and fix some index $i$. 
The following result is the main tool in checking that specific 
manifolds are not w.g.s.c.

\begin{proposition}\label{T:wgsc}
Any exhaustion as above of the w.g.s.c. manifold $W^n$ has a
refinement for which consecutive terms fulfill the conditions: 
\begin{enumerate}
\item all pairs $(T_{i+1}^n,T_i^n)$ are stably-compressible. 
\item if $S_{i,j}^{n-1}$ denote the components of $\del T_i^n$ the map 
$*\pi_1(S_{i,j}^{n-1})\to \pi_1(T_i^n)$ induced by the inclusion 
is surjective. 
\end{enumerate} 
\end{proposition}

\begin{proof} 
 Since $W^n$ is w.g.s.c. there 
exists a compact 1-connected submanifold $M^n$ of $W^n$ such 
that $T_i^n\subset M^n$. We can suppose $M^n\subset T_{i+1}^n$, without 
loss of generality. From now on we will focus on the pair 
$(T_{i+1}^n, T_i^n)$  and suppress the index $i$, and denote it 
$(T', T)$, for the sake of notational  simplicity.

\begin{lemma}\label{stably}
The pair $(T', T)$ is stably compressible.
\end{lemma}

\begin{proof}
Let 
$\varphi:\pi_1(\partial T)\to \pi_1(T)$ and  
$\psi:\pi_1(\partial T)\to \pi_1(T'-{\rm int}(T))$
be the homomorphisms induced by the inclusions 
$\partial T\hookrightarrow T$, 
and $\partial T\hookrightarrow T'-{\rm int}(T)$. If $\partial T$ has several
components then we choose base points in each component and 
set $\pi_1(\partial T)=*_j\pi_1(S_j)$ for notational simplicity. 

\vspace{0.2cm}
\noindent 
Let us consider a handlebody decomposition of $M^n-{\rm int}(T)$ (respectively
a connected component) from $\partial T$,
$$M^n-{\rm int}(T)= \left(\partial T \times[0,1]\right) 
\bigcup_{\lambda=1}^{n-1} \left(\cup_j 
h_j^{\lambda}\right), $$
where $h_j^{\lambda}$ is a handle of index
$\lambda$. We suppose the handles are attached in increasing order of
their index.  Since the distinct components of $\partial T$ are not
connected outside $T$ (by remark 3.13) the 1-handles which are added have the
extremities in the same connected component of $\partial T$.  Set
$M_2^n\subset M^n$ (respectively $M_1^n$) for the submanifold obtained
by attaching to $T$ only the handles $h_j^{\lambda}$ of index
$\lambda\leq 2$ (respectively those of index $\lambda\leq 1$).  Then
$\pi_1(M_2^n)=0$, because adding higher index handles does not affect
the fundamental group and we know that $\pi_1(M^n)=0$.

\begin{lemma}
The pair $(T', M_1^n)$ is strongly compressible.
\end{lemma}

\begin{proof}
Let $\{\gamma_j\}_{j=1,p}\subset \partial M_1^n$ be the set of 
attaching circles for the 2-handles of $M_2^n$ and $\{\delta^2_j\}_{j=1,p}$ 
be the corresponding core of the 2-handles $h^2_j$ ($j=1,p$). 
Since  $\delta^2_j$ is a  
2-disk embedded in $M^n-{\rm int}(T)$  it follows that 
the homotopy class $[\gamma_j]$ vanishes in $\pi_1(M^n-{\rm int}(T))$. 

\vspace{0.2cm}
\noindent 
Let $\Gamma\subset \pi_1(\partial M_1^n)$ be the normal subgroup 
generated by the homotopy classes of the curves $\{\gamma_j\}_{j=1,p}$
which are contained in $\partial M_1^n$. Notice that this amounts to 
picking base points which are joined to the loops.  Therefore 
the image of $\Gamma$ under the map 
$\pi_1(\partial M_1^n)\to \pi_1(M_2^n-{\rm int}(M_1^n))$, 
induced by the inclusion, is zero. In particular its image in 
$\pi_1(T'-{\rm int}(M_1^n))$ is zero.  

\vspace{0.2cm}
\noindent 
On the other hand the images of the 
classes $[\gamma_j]$ in $\pi_1(M_1^n)$  
normally generate all of the group $\pi_1(M_1^n)$ because 
$\pi_1(M_2^n)=1$.  

\vspace{0.2cm}
\noindent 
These two properties are equivalent to the strong compressibility 
of the pair  $(M_2^n, M_1^n)$ which in turn implies that of 
$(T', M_1^n)$.
\end{proof}

\noindent  
{\em Rest of the proof of Lemma~\ref{stably}:}  Assume now that the number of 1-handles $h_j^1$ is $r$. 
Notice that $\pi_1(M_1^n-{\rm int}(T))\cong \pi_1(\partial T)*{\bf F}(r)$, 
because it can be obtained from $\partial T\times [0,1]$ by adding 
1-handles and the 1-handles we added do not join distinct 
boundary components, so that each one contributes with a free factor. 
In particular the inclusion 
$\partial T\hookrightarrow M_1^n-{\rm int}(T)$ induces a 
monomorphism 
$\pi_1(\partial T)\hookrightarrow \pi_1(\partial T)*{\bf F}(r)$. 
Observe that $M_1^n-{\rm int}(T)$ can also be obtained from 
$\partial M_1^n \times[0,1]$ by adding $(n-1)$-handles hence 
the inclusion $\partial M_1^n \hookrightarrow M_1^n-{\rm int}(T)$
induces the isomorphism 
$\pi_1(\partial M_1^n)\cong \pi_1(\partial T)*{\bf F}(r)$. 
The same reasoning gives the isomorphism 
$\pi_1(M_1^n)\cong \pi_1(T)*{\bf F}(r)$. 

\vspace{0.2cm}
\noindent 
In particular we can view the subgroup 
$\Gamma$ as a subgroup of $\pi_1(\partial T)*{\bf F}(r)$.
The previous lemma tells us that $\Gamma$ lies in the kernel of 
$\pi_1(\partial T)*{\bf F}(r)\to \pi_1(T'-{\rm int}(T))$ 
and also projects epimorphically onto 
$\pi_1(T)*{\bf F}(r)$. The identification of the respective maps 
with the morphisms induced by inclusions yields our claim.
\end{proof}
This finishes the proof of Proposition \ref{T:wgsc}.
\end{proof}

\noindent 
Conversely, assume that $W^n$ has an exhaustion in which consecutive 
pairs are stably compressible. Then it is sufficient to show the 
following:

\begin{proposition}\label{weak}
If $(T',T)$ is a stably compressible pair of $n$-manifolds and 
$n\geq 5$ then 
$T\subset {\rm int}(M^n)\subset T'$ where 
$M^n$ is a compact submanifold with $\pi_1(M^n)=0$. 
\end{proposition}
\begin{proof}

\noindent  
One can  realize the homomorphism $\beta:{\bf F}(r)\to \pi_1(T'-{\rm int}(T))$ 
by a disjoint union of  bouquets of circles $\vee_{}^r S^1\to T'-{\rm int}(T)$. 
There is one bouquet in each connected component of $T'-{\rm int}(T)$.
One joins each wedge point to the unique
connected component of $\partial T$ for which that is possible by an arc, 
and set $M_1^n$ for the manifold obtained from $T$ by adding 
a regular neighborhood of the bouquets $\vee_{}^r S^1$ in $T'$ 
(plus the extra  arcs). This is equivalent to adding 1-handles 
with the induced framing.

\begin{lemma}
The kernel $\ker \psi*\beta\subset \pi_1(\partial M_1^n)$ is 
normally generated by a finite number of elements 
$\gamma_1,\gamma_2,..., \gamma_p$. 
\end{lemma}

\begin{proof}
Consider a finite presentation ${\bf F}(k)/H\to \pi_1(\partial T)$. 
We know that  $ \pi_1(\partial M_1^n)= \pi_1(\partial T)*{\bf F}(r)$. 
Furthermore the composition map 
\[ \lambda: {\bf F}(k)* {\bf F}(r)\to \pi_1(\partial T)*{\bf F}(r)
\stackrel{\psi*\beta}{\to} \Gamma, \]
is surjective (since $\beta$ is). The first map is the free 
product of the natural projection with the identity. 
Therefore  ${\bf F}(k+r)/\ker \lambda\cong \Gamma$ is a 
presentation of the group $\Gamma$. 
The Theorem of Neumann (see \cite{Ba}, p.52) states that any
presentation on finitely many generators of a finitely presented group
has a presentation on these generators  with only finitely many 
of the given relations. Applying this to  $\Gamma$ 
one derives that there exist finitely many elements which 
normally generate $\ker \lambda$ in ${\bf F}(k+r)$. 
Then the images of these elements  in 
$\pi_1(\partial T)*{\bf F}(r)$ normally generate $\ker \psi*\beta$
(the projection ${\bf F}(k)* {\bf F}(r)\to \pi_1(\partial T)*{\bf
  F}(r)$ is surjective). This yields the claim. 
\end{proof}

\begin{lemma}
The elements $\gamma_i$  are also in the kernel of 
$\pi_1(\partial M_1^n)\to \pi_1(T'-{\rm int}(M_1^n))$.
\end{lemma}

\begin{proof}
The map 
$\pi_1(T'-{\rm int}(M_1^n))\to\pi_1(T'-{\rm int}(T))$ induced by the inclusion 
is injective  because $T'-{\rm int}(T)$   
is obtained from $T'-{\rm int}(M_1^n)$ by adding $(n-1)$-handles (dual to 
the 1-handles from which one gets $M_1^n$ starting from $T$), and 
$n\geq 5$. 
Thus the map $\pi_1(\partial M_1^n)\to \pi_1(T'-{\rm int}(T))$ factors
through 
\[ \pi_1(\partial M_1^n)\to \pi_1(T'-{\rm int}(M_1^n))\to\pi_1(T'-{\rm int}(T)),
\]
and any element in the kernel must be in the kernel of the first map,
as stated.
\end{proof}

\noindent 
The dimension restriction $n\geq 5$ implies that
we can assume $\gamma_j$ are represented by embedded 
loops having only the base point in common.
Then $\gamma_j$ bound singular 2-disks $D^2_j\subset T'-{\rm int}(M_1^n)$.
By a general position argument, one can 
arrange such that the  2-disks $D^2_j$  are embedded in
$T'-{\rm int}(M_1^n)$ 
and have disjoint interiors.
 
\vspace{0.2cm}
\noindent 
As a consequence the manifold $M^n$ obtained from $M_1^n$ by attaching 
2-handles along the $\gamma_j$'s (with the induced framing) 
can be embedded in $T'-{\rm int}(T)$. Moreover $M^n$  is 
a compact manifold whose fundamental group 
is the quotient of  $\pi_1(M_1^n)$ by the subgroup normally generated 
by the elements $\varphi*1(\gamma_j)$'s. 
The group $\varphi *1(\ker \psi*\beta)$ is normally generated 
by the elements $\varphi*1(\gamma_j)$. By hypothesis 
the pair $(\varphi *1,\psi*\beta)$ is compressible hence 
$\varphi *1(\ker \psi*\beta)$ contains $\varphi*1(\pi_1(\partial
M_1^n))=  \varphi*1(\pi_1(\partial T)*{\bf F}(r))=
\varphi(\pi_1(\partial T))*{\bf F}(r)$. 
Next 
\[\frac{\pi_1(M_1^n)}{\varphi *1(\ker \psi*\beta)}=
\frac{\pi_1(T)*{\bf F}(r)}{\varphi(\pi_1(\partial T))*{\bf F}(r)}\cong
\frac{\pi_1(T)}{\varphi(\pi_1(\partial T))}=1,\]
since $\varphi$ has been supposed surjective. 
Therefore the quotient of $\pi_1(M_1^n)$ by the subgroup normally
generated by the elements $\varphi*1(\gamma_j)$ is trivial.
\end{proof}

\subsection{End 1-compressibility is trivial for $n\geq 5$}
We defined an infinite sequence of obstructions (namely
$k$-compressibility for each $k$) to the w.g.s.c. However the first 
obstruction is trivial i.e. equivalent to the simple connectivity, 
in dimensions $n\neq 4$. 
In fact the main result of this section establishes the following:  

\begin{proposition}\label{comprestrivial}
End 1-compressibility and simple connectivity (s.c.) are equivalent 
for open $n$-manifolds of dimension $n\geq 5$. 
\end{proposition}

\begin{proof}
We first consider a simpler case: 
\begin{proposition}
The result holds in the case of a manifold $W^n$ of dimension at least
$5$ with one end. 
\end{proposition}

\begin{proof} In this case  $W^n$ has an exhaustion $T_i$ with $\partial
T_i$ connected for all $i$.

\begin{lemma}\label{T:expand} $W^n$ has an exhaustion such that the map
$\varphi:\pi_1(\del T_i)\to \pi_1(T_i)$ induced by inclusion is a
surjection for each $i$.
\end{lemma}

\begin{proof} As $W^n$ is simply connected, by taking a refinement we can
assume that each inclusion map $\pi_1(T_i)\to\pi_1(T_{i+1})$ is the
zero map. As usual we denote $T_i$ and $T_{i+1}$ by $T$ and $T'$
respectively.

\vspace{0.2cm}
\noindent 
Now, take a handle-decomposition of $T$ starting with the boundary
$\del T$. Suppose the core of each $1$-handle of this decomposition is
homotopically trivial in $T$, then it is immediate that $\varphi$ is a
surjection. We will enlarge $T$ by adding some $1$-handles and
$2$-handles (that are embedded in $T'$) at $\del T$ in order 
to achieve this.

\vspace{0.2cm}
\noindent 
Namely, let $\gamma$ be the core of a $1$-handle. By hypothesis, there
is a disc $D^2$ in $T'$ bounding the core of each of the $1$-handles,
which we take to be transversal to $\del T$. As the dimension of $W^n$
is at least $5$, $D^2$ can be taken to be embedded. Notice that the
2-disks corresponding to all 1-handles can also be made disjoint, by 
general position.  Thus $D^2$ intersects
$T'-{\rm int}(T)$ in a collection of embedded disjoint planar surfaces. The
neighborhood of each disc component of this intersection can be
regarded as a $2$-handle (embedded in $T'$) which we add to $T$ at
$\del T$. For components of $D^2\cap (T'-{\rm int}(T))=D^2-{\rm int}(T)$ with more
than one boundary component, we take embedded arcs joining distinct
boundary components. We add to $T$ a neighborhood of each arc, which
can be regarded as a $1$-handle. After doing this for a finite
collection of arcs, $D^2-D^2\cap T$ becomes a union of discs. Now we add
$2$-handles as before. The disc $D^2$ that $\gamma$ bounds is now in
$T$.

\vspace{0.2cm}
\noindent 
Further, the dual handles to the handles added are of dimension at
least $3$. In particular we can extend the previous
handle-decomposition to a new one for (the new) $T$ starting at (the
new) $\del T$ with no new $1$-handles. Thus, after performing the
above operation for the core of each $1$-handle of the original handle
decomposition, the core of each $1$-handle of the resulting
handle-decomposition of $T$ starting at $\del T$ bounds a disc in
$T$. Thus $\phi$ is a surjection.  \end{proof}

\vspace{0.2cm}
\noindent 
The above exhaustion is in fact $1$-compressible by the following
algebraic lemma.

\begin{lemma}\label{T:VK} Suppose that we have a square of maps
verifying the Van Kampen theorem:

$$\begin{CD}
A               @>\psi>>        C \\
@V{\varphi}VV           @VV{\gamma}V\\
B               @>\beta>>       D
\end{CD}$$

\noindent Let $\xi: A\to D$.  Suppose also that $\varphi$ is surjective and $\beta$
is the zero map. Then $\psi(A)\subset {\cal N}(\psi(\ker \varphi))=
{\cal N}(\psi(\ker \varphi), C)$.
\end{lemma}
\begin{proof} Observe that $\xi(\ker \varphi)=0$, hence $\psi(\ker
  \varphi)\subset \ker \gamma$.
Hence we can define another diagram with $A'=A/\ker \varphi$,
$B'=B$, $C'=C/{\cal N}(\psi(\ker \varphi))$, $D'=D$:
 
$$\begin{CD}
A'=A/\ker \varphi        @>\psi>> C'=C/{\cal N}(\psi(\ker \varphi))\\
@V{\varphi}VV                   @VV{\gamma}V\\
B'=B            @>\beta>>       D'=D
\end{CD}$$

\vspace{0.2cm}
\noindent  
Notice that the map $A/\ker \varphi \to C/{\cal N}(\psi(\ker \varphi))$ is
well-defined.  Again this diagram verifies the Van Kampen theorem.
For this diagram, the induced $\varphi$ is an isomorphism.  It is
immediate then that the universal (freest) $D'$ must be $C'$.  In fact
$D'=C'*A'/ {\cal N}( \{\psi(a)a^{-1}, a\in A'\})$.  Consider the map
$C'\to C'*A'/{\cal N} ( \{\psi(a)a^{-1}, a\in A'\})\to C'$, where the second
arrow consists in replacing any occurrence of $a$ by the element
$\psi(a)$ and taking the product in $C'$. This composition is the
identity and the first map is a surjection, 
hence the map $C'\to D'$ is an isomorphism.

\vspace{0.2cm}
\noindent 
The map induced by $\xi$ is zero since $\beta$ is the zero.  But the map
$\xi :A'\to D'$ is the map $\psi:A'\to C'$ followed by an isomorphism,
hence $\psi(A')=0$. This is equivalent to $\psi(A)\subset
{\cal N}(\psi(\ker\varphi))$.
\end{proof}

\noindent 
Note that the above lemma is purely algebraic, and in particular
independent of dimension. The two lemmas immediately give us the
proposition for one-ended manifolds $W^n$ with $n\geq 5$. 
\end{proof}

\noindent {\em The general case.}
We now consider the general case of a simply-connected open manifold
$W^n$ of dimension at least $5$, with possibly more than one end. We
shall choose the exhaustion $T_i$ with more care in this case. 

\vspace{0.2cm}
\noindent 
We will make use of the following construction several times.  Start
with a compact submanifold $A^n$ of codimension $0$, with possibly more
than one boundary component. Assume for simplicity (by enlarging $A^n$ if
necessary) that no complementary component of $A^n$ is pre-compact.  As
$W^n$ is simply connected, we can find a compact submanifold $B^n$
containing $A^n$ in its interior such that the inclusion map on
fundamental groups is the zero map. Further, we can do this by
thickening and then adding the neighborhood of a $2$-complex, i.e., a
collection of $1$-handles and $2$-handles. Namely, for each generator
$\gamma$ of $\pi_1(A^n)$, we can find a disc $D^2$ that $\gamma$ bounds,
and then add $1$-handles and $2$-handles as in
lemma~\ref{T:expand}. Thus, as $n\geq 5$, the boundary components
of $B^n$ correspond to those of $A^n$. We repeat this with $B^n$ in place of
$A^n$ to get another submanifold $C^n$.

\vspace{0.2cm}
\noindent 
Observe that as a consequence of this and the simple-connectivity of
$W^n$, the inclusion map $\pi_1(A^n\cup V^n)\to \pi_1(B^n\cup V^n)$ is the zero
map for any component $V^n$ of $W^n - A^n$. More generally if $Z^n\subset V^n$,
then $\ker(\pi_1(A^n\cup Z^n)\to \pi_1(B^n\cup Z^n)) = \ker(\pi_1(A^n\cup Z^n)\to
\pi_1(B^n\cup Z^n\cup (W^n - V^n)))$. Similar results hold with $B^n$ and $C^n$ in
place of $A^n$ and $B^n$.

\vspace{0.2cm}
\noindent 
Now start with some $A_1^n$ as above and construct $B_1^n$ and
$C_1^n$. Thicken $C_1^n$ slightly to get $T=T_1$. We will eventually
choose a $T_2=T'$, but for now we merely note that it can (and so  it will)
be chosen in such a manner that the inclusion map $\pi_1(T)\to
\pi_1(T')$ is the zero map. Let $S_j,j=1,\dots,n$ be the boundary
components of $T$. Let $X_j$ be the union of the component of $T -
{\rm int}(B^n)$ containing $S_j$ and $B^n$, and define $Y_j$ analogously with
$C^n$ in place of $B^n$. Denote the image of $\pi_1(X_j)$ in $\pi_1(Y_j)$
by $\overline{\pi}(X_j)$. We then have a natural map
$\varphi:\pi_1(S_j)\to\overline{\pi}(X_j)$. Let $V_j$ be the component of
$T'-{\rm int}(T)$ containing $S_j$.

\begin{lemma}\label{T:exp2} By adding $1$-handles and $2$-handles to
$S_j$, we can ensure that $\pi_1(S_j)$ surjects onto $\overline{\pi}(X_j)$.
\end{lemma}
\begin{proof} The proof is essentially the same as that of
lemma~\ref{T:expand}. We start with a handle-decomposition for $X_j$
starting from $S_j$. We shall ensure that the image in $\overline{\pi}(X_j)$
of the core of each $1$-handles is trivial. Namely, for each core, we
take a disc $D$ that it bounds.  By the above remarks, we can, and do,
ensure that the disc lies in $Y_j$, and in particular does not
intersect any boundary component of $T$ except $S_j$. As in
lemma~\ref{T:expand} we may now add $1$-handles and $2$-handles to
$S_j$ to achieve the desired result.  \end{proof}

\noindent 
Notice that the changes made to $T$ in the above lemma do not affect
$S_k$, $X_k$ and $Y_k$ for $k\neq j$. Hence, by repeated application
of the above lemma, we can ensure that all the maps
$\pi_1(S_j)\to\overline{\pi}(X_j)$ are surjections. Also notice that the
preceding remarks show that
$\ker \varphi=\ker(\pi_1(S_j)\to\pi_1(T))$. Now take $A=\pi_1(S_j)$,
$B=\overline{\pi}(X_j)$ and $C=\pi_1(V_j)$, and let $D$ be the image of
$\pi_1(V_j\cup X_j)$ in $\pi_1(V_j\cup Y_j)$. Then, by the preceding
remarks and lemma~\ref{T:exp2}, the diagram

$$\begin{CD}
A               @>\psi>>        C \\
@V{\varphi}VV           @VV{\gamma}V\\
B               @>\beta>>       D
\end{CD}$$
\noindent satisfies the hypothesis of the lemma~\ref{T:VK}. The 1-compressibility
for the pair $(T',T)$ follows.

\vspace{0.2cm}
\noindent 
Now we continue the process inductively. Suppose $T_k$ has been
defined, choose $A_{k+1}$ so that it contains $T_k$ and also in such a
manner as to ensure that $A_i$'s exhaust $M$. Then find $B_{k+1}$,
$C_{k+1}$ and $T_{k+1}$ as above. The rest follows as above.
\end{proof}
\begin{remark} For the case of simply-connected, one-ended (hence contractible)
$3$-manifolds, a theorem of Luft says that $M$ can be exhausted by a
union of homotopy handlebodies. These satisfy the conclusion of
lemma~\ref{T:expand}, hence the proposition still holds. More
generally, we can apply the sphere theorem to deduce that we have an
exhaustion by connected sums of homotopy handlebodies. It follows that
each pair $(T,T')$ of this exhaustion is $1$-compressible as we can
decompose $T$ and consider each component separately without affecting
$\varphi$ or $\psi$.
\end{remark}

\section{Examples of contractible manifolds}
\subsection{Uncountably many Whitehead-type manifolds}

\noindent Recall the following definition from \cite{Wr}:

\begin{definition}
A {\em Whitehead link} $T^n_0\subset T^n_1$ is a 
null-homotopic embedding of the solid torus 
$T^n_0$ in the (interior of the) 
unknotted solid torus $T^n_1$ lying in $S^n$ such that 
the pair $(T_1^n,T_0^n)$ is (boundary) incompressible. 
\end{definition}

\noindent 
The solid $n$-torus is $T^n=D^2\times S^1\times S^1...\times S^1$. 
By iterating the ambient homeomorphism which sends  $T^n_0$ onto $T^n_1$ 
one obtains an ascending sequence $T^n_0\subset T^n_1\subset T^n_2\subset ...$ 
whose union is called a Whitehead-type $n$-manifold. 
A Whitehead-type manifold is open contractible and not s.c.i. 
D.G.~Wright (\cite{Wr}) gave a recurrent procedure 
to construct many Whitehead links in dimensions $n\geq 3$. 
One shows below  that this construction provides uncountably many 
distinct contractible manifolds.

\vspace{0.2cm}
\noindent 
We introduce an invariant for pairs of solid tori which generalizes the 
wrapping number in dimension 3. Moreover this provides invariants 
for open manifolds of Whitehead-type answering a question raised 
in \cite{Wr}.

\begin{definition}
A {\em spine} of the solid torus $T^n$ is an embedded 
$t^{n-2}=\{*\}\times S^1\times S^1\times...\times S^1\subset T^n$ 
having a trivial normal bundle in $T^n$. This gives $T^n$ the structure 
of a trivial  2-disk bundle over $t^{n-2}$. 
\end{definition}

\begin{remark}
Although the spine is not uniquely defined, its isotopy class 
within the solid torus is. 
\end{remark}

\noindent 
Consider a pair of solid tori $T_0^n\subset T^n$. We fix some spine 
$t^{n-2}$ for $T^n$. To specify the embedding of $T^n_0$ is the same 
as giving the embedding of a spine $t^{n-2}_0$ of $T^{n}_0$ in $T^n$. 
The isotopy class of the embedding $t^{n-2}_0\hookrightarrow T^n$ is 
therefore uniquely defined by the pair.
Let us pick-up a Riemannian metric $g$ on the 
torus $T^n$ such that $T^{n}$ is identified with the regular 
neighborhood of radius $r$
around $t^{n-2}$. We denote this by $T^{n}=t^{n-2}[r]$, and 
suppose for simplicity that $r=1$. Then $t^{n-2}[\lambda]$ for 
$\lambda\leq 1$ will denote the radius $\lambda$ tube around $t^{n-2}$
in this metric. 
\begin{definition}
The {\em wrapping number} of the Whitehead link $T^n_0\subset T^n$ is 
defined as follows:
\[ w(T^n,T^n_0)=\lim_{\varepsilon\to 0} \;\;\; \inf_{t^n_0\in {\cal I}(
t^n_0\subset t^{n-2}[\varepsilon])}\;\;\;
\frac{{\rm vol}(t^{n-2}_0)}{{\rm vol}(t^{n-2})}, \]
where ${\cal I}(t^n_0\subset T^n)$ is the set of all embeddings 
of the spine $t^{n-2}_0$ of $T^n_0$ in the given  isotopy class, and ${\rm vol}$ 
is the $(n-2)$-dimensional volume.
\end{definition}
\begin{remark}
Notice that a priori this definition might depend on the particular choice 
of the spine $t^{n-2}$ and on the metric $g$. 
\end{remark}
\begin{proposition}
The wrapping number is a topological invariant of the pair 
$(T^n,T^n_0)$. 
\end{proposition}
\begin{proof}
There is a natural projection map on the spine $\pi:T^n\to t^{n-2}$, which 
is the fiber bundle projection of $T^n$ (with fiber a 2-disk). 
When both $T^n$ and $t^{n-2}$ are 
fixed then such a projection map is also defined only up to  isotopy. 
Set therefore 
\[l(T^{n},T^{n}_0)=\inf_{t^n_0\in {\cal I}(t^n_0\subset T^n)} \;\;
\inf_{x\in t^{n-2}} \sharp\left\{\pi^{-1}(x)\cap t^{n-2}_0 \right\}.\]
Since $\inf_{x\in t^{n-2}}\sharp\left\{\pi^{-1}(x)\cap t^{n-2}_0 \right\}$ 
does not depend on the particular projection map (in the fixed isotopy class)
this number  represent a topological invariant of the pair $(T^n,T^n_0)$.
Hence  the claim follows from the following result:
\begin{proposition}
$w(T^n,T^n_0)=l(T^{n}, T^{n}_0).$
\end{proposition}
\begin{proof}
Consider a position of $t^{n-2}_0$ for which the minimum 
value $l(T^{n},T^{n}_0)$ is attained. 
A small isotopy make  $t^{n-2}_0$ transversal to $\pi$. 
Then, for this precise position of $t^{n-2}_0$ there exists some 
number $M$ such that 
\[  \sharp\left\{\pi^{-1}(x)\cap t^{n-2}_0 \right\} \leq M, \mbox{ for any
} x\in t^{n-2}. \]
Denote by $\mu$ the Lebesgue measure on $t^{n-2}$. 
\begin{lemma}\label{good}
For any $\varepsilon > 0$ one can move $t^{n-2}_0$ in $T^n$ by an ambient 
isotopy such that the following conditions are fulfilled:
\[   \sharp\left\{\pi^{-1}(x)\cap t^{n-2}_0 \right\} \leq M, \mbox{ for any
} x\in t^{n-2}.\]
\[ \mu\left(\left\{ x\in t^{n-2}|\;\; \sharp\left\{\pi^{-1}(x)\cap t^{n-2}_0 
\right\} > l(T^{n},T^{n}_0)\right\}\right) < \varepsilon.\]
\end{lemma}
\begin{proof}
The set 
$U=\left\{ x\in t^{n-2}|\;\; \sharp\left\{\pi^{-1}(x)\cap t^{n-2}_0 \right\} = l(T^{n},T^{n}_0)\right\}$ is an open subset of positive measure. 
Consider then a flow $\varphi_t$ on the torus $t^{n-2}$ which expands 
a small ball contained in $U$ into the complement of a measure 
$\varepsilon$ set (e.g. a  small tubular neighborhood of a spine of the 
1-holed torus). Extend this flow as $1_{D^2}\times \varphi_t$ all over 
$T^n$ and consider 
its action on $t^{n-2}_0$.
\end{proof}

\begin{lemma}
$w(T^n,T^n_0)\leq l(T^{n}, T^{n}_0).$
\end{lemma}
\begin{proof}
The map $\pi$ is the projection of the metric tube around $t^{n-2}$ on 
its spine, hence the Jacobian $Jac(\pi|_{t^{n-2}_0})$ has bounded norm 
$|Jac(\pi|_{t^{n-2}_0})| \leq 1$. 
It follows that 
\[ \frac{{\rm vol}(t^{n-2}_0)}{{\rm vol}(t^{n-2})} = 
\frac{\int \pi^* d \mu}{\int d\mu}
\leq \frac{\int_{\pi^{-1}(U)} |Jac(\pi|_{t^{n-2}_0})| d \mu}{\int_U
  d\mu} + M\varepsilon \leq \]
\[ \leq l(T^{n}, T^{n}_0)(1-\varepsilon) + M\varepsilon,\]
for any $\varepsilon >0$, hence the claim.
\end{proof}

\begin{lemma}
$w(T^n,T^n_0)\geq l(T^{n}, T^{n}_0).$
\end{lemma}

\begin{proof}
Set $\lambda_t:t^{n-2}[\delta]\to t^{n-2}[t\delta]$ for 
the map given in coordinates by  $\lambda_t(p,x)=(tp,x)$, $p\in D^2, 
x\in t^{n-2}$.  Here the projection $\pi$ provides a global trivialisation 
of $t^{n-2}[\delta]\subset T^n$. Then 
\[\lim_{t\to 0} |Jac(\pi|_{t^{n-2}_0}\circ \lambda_t)|=1.\]
Therefore for $t$ close enough to 0 one derives  
\[\lim_{t\to 0}\frac{{\rm vol}(\lambda_t(t^{n-2}_0))}{{\rm vol}(t^{n-2})}=
\lim_{t\to 0} \frac{\int_{\lambda_t(t^{n-2}_0)} |Jac(\pi|_{t^{n-2}_0})| d \mu}{\int_
{t^{n-2}} d\mu}\geq l(T^{n}, T^{n}_0).\]
Since the position of $t^{n-2}_0$ was chosen arbitrary, this inequality 
survives after passing to the infimum and the claim follows.
\end{proof}
\end{proof}\end{proof}
\begin{theorem}\label{uncountable}
There exist uncountably many Whitehead-type manifolds for $n\geq 5$. 
\end{theorem}
\begin{proof} 
The proof here follows the same pattern as that given by 
McMillan (\cite{McM}) for the 3-dimensional case. 
Let us establish first the following useful property of the wrapping number:
\begin{proposition}
If $T_0^n\subset T_1^n\subset T_2^n$ then 
$w(T^n_2,T^n_0)=w(T^n_2,T^n_1)w(T^n_1,T^n_0)$.
\end{proposition}
\begin{proof}
This is a consequence of the two lemmas below:
\begin{lemma}
$l(T^n_2,T^n_0)\leq l(T^n_2,T^n_1)l(T^n_1,T^n_0)$.
\end{lemma}
\begin{proof}
Consider  $t_0^{n-2}\subset T^n_1\subset T^n_2$, where $T_1^n$ is a 
very thin tube around $t_1^{n-2}$, and the two projections 
to $\pi_2:t^{n-2}_1\to t^{n-2}_2$ and 
$\pi_1:t^{n-2}_0\to t^{n-2}_1$ respectively.
Using Lemma \ref{good} one can assume that the conditions
\[ \mu\left(\left\{ x\in t_1^{n-2}|\;\; \sharp\left\{\pi_1^{-1}(x)\right\} 
= l(T^{n}_1,T^{n}_0)\right\}\right) > 1 - \varepsilon,\]
\[ \mu\left(\left\{ x\in t_2^{n-2}|\;\; \sharp\left\{\pi_2^{-1}(x)\right\} = l(T^{n}_2,T^{n}_1)\right\}\right) > 1 - \varepsilon,\]
hold. For small enough $\epsilon$ one  derives that 
\[ \mu\left(\left\{ x\in t_2^{n-2}|\;\; \sharp\left\{(\pi_2\circ\pi_1)^{-1}(x)\right\} = l(T^{n}_2,T^{n}_1)l(T^{n}_1,T^{n}_0)\right\}\right) > 0.\]
This proves that the minimal cardinal of the 
$(\pi_2\circ\pi_1)^{-1}(x)$ is not greater than $l(T^{n}_2,T^{n}_1)l(T^{n}_1,T^{n}_0)$, hence the claim.
\end{proof}
\begin{lemma}
$w(T^n_2,T^n_0)\geq w(T^n_2,T^n_1)w(T^n_1,T^n_0)$.
\end{lemma}
\begin{proof}
We can assume that  $w(T^n_2,T^n_1)\neq 0$.
Consider  an embedding of the $(n-2)$-torus $s^{n-2}_1\subset T_2=t^{n-2}_2[\varepsilon]$  for which the  value 
of $\frac{{\rm vol}(t^{n-2}_1)}{{\rm vol}(t^{n-2}_2)}$ (as 
function of $t^{n-2}_1$) is closed to the 
infimum in the isotopy class. We will assume that in all formulas
below the tori lay in their respective isotopy classes.
Then 
\[ \left(\inf_{t_0^{n-2}\subset s_1^{n-2}[2\varepsilon]}
\frac{{\rm vol}(t^{n-2}_0)}{{\rm vol}(s^{n-2}_1)}\right)\left( \inf_{t_1^{n-2}\subset t_2^{n-2}[\varepsilon]} \frac{{\rm vol}(t^{n-2}_1)}{{\rm vol}(t^{n-2}_2)}\right) \leq \]
\[\leq \left(\inf_{t_0^{n-2}\subset s_1^{n-2}[2\varepsilon]}
\frac{{\rm vol}(t^{n-2}_0)}{{\rm vol}(t^{n-2}_1)}\right) 
\frac{{\rm vol}(s^{n-2}_1)}{{\rm vol}(t^{n-2}_2)} = \]
\[= \inf_{t_0^{n-2}\subset s_1^{n-2}[2\varepsilon]}
\frac{{\rm vol}(t^{n-2}_0)}{{\rm vol}(t^{n-2}_2)} \leq 
\inf_{t_0^{n-2}\subset t_2^{n-2}[\varepsilon]}
\frac{{\rm vol}(t^{n-2}_0)}{{\rm vol}(t^{n-2}_2)}. \]
The last inequality follows from the fact that 
$s_1^{n-2}[2\varepsilon]\subset t_2^{n-2}[\varepsilon]$.
In fact  $w(T^n_2,T^n_1)\neq 0$ implies that 
$s_1^{n-2}$ intersects any 2-disk $D^2\times \{*\}$ (i.e. 
any fiber of the projection $\pi_2:T^n_2\to t^{n-2}_2$) of 
$T^n_2=t^{n-2}_2[\varepsilon]$ in at least 
one point. Then the transversal  disk $D^2\times \{*\}$ of radius 
$\varepsilon$ is  therefore contained  in 
the tube $s_1^{n-2}[2\varepsilon]$ of 
radius $2\varepsilon$ around $s_1^{n-2}$, establishing the 
claimed inclusion.  

\vspace{0.2cm}
\noindent 
On the other hand the following holds
\[ \lim_{\varepsilon\to 0}\;\;\;\inf_{t_0^{n-2}\subset s_1^{n-2}[2\varepsilon]}\;\;
\frac{{\rm vol}(t^{n-2}_0)}{{\rm vol}(s^{n-2}_1)}=
w(T^n_1,T^n_0), \]
due to the topological invariance of the wrapping number. 
Letting $\varepsilon$ go to 0 in the previous inequality
yields the claim.  
\end{proof}
\end{proof}
\begin{proposition}\label{wrconst}
There exist Whitehead links whose wrapping number has the form 
$2^{n-2}p$ for any natural number $p$.
\end{proposition}
\begin{proof}
The claim is well-known for $n=3$.
One uses Wright's construction (\cite{Wr}) of Whitehead links 
by induction on the 
dimension. If $T_0^n\subset T_1^n$ is a Whitehead link then 
set $T^{n+1}=T_1^n\times S^1$. Consider the projection $q$ of the 
solid torus $T_0^n\times S^1\cong D^2\times S^1\times...\times S^1$ 
onto $D^2\times S^1$ (the first and the last factors). Choose some 
Whitehead link  $L^3\subset D^2\times S^1$, and set then 
$Q^{n+1}=q^{-1}(L^3)$. The pair $Q^{n+1}\subset T^{n+1}$ 
is a Whitehead link of dimension $n+1$. 
The Proposition then is an immediate consequence of:
\begin{lemma}
$w(T^{n+1},Q^{n+1})=w(T_1^{n},T_0^{n})w(D^2\times S^1,L^3)$.
\end{lemma}
\begin{proof}
From the multiplicativity of $w$ and the triviality of the projection $q$
it is sufficient to prove that 
$w(T_1^{n}\times S^1,T_0^{n}\times S^1)=w(T_1^{n},T_0^{n})$. 
This formula can be checked directly using $l$ instead of $w$. 
\end{proof}\end{proof}
\begin{proposition}
For any sequence ${\bf p}=p_0,p_1,...$ of positive integers 
consider a Whitehead-type manifold $W^n({\bf p})=\cup_{k=1}^{\infty}T^n_k$, 
where $w(T^n_{k+1}, T^n_k)=2^{n-2}p_k$. 
If the sequences ${\bf p}$ and ${\bf q}$ have infinitely many 
non-overlapping prime factors then the manifolds 
$W^n({\bf p})$ and $W^n({\bf q})$ are not PL homeomorphic. 
\end{proposition}
\begin{proof}
The proof is similar to that of (\cite{McM}, p.375).
Set  $W^n({\bf p})=\cup_{k=1}^{\infty}T^n_k$, 
$W^n({\bf q})=\cup_{k=1}^{\infty}s^n_k$, where 
$T^n_k$,$\widetilde{T}^n_k$, are tori, as above. 
If $h:W^n({\bf q})\to W^n({\bf p})$ is a PL homeomorphism, 
there exist integers $j,k$ such that 
$T^n_0\subset int(h(\widetilde{T}^n_j))$, $q_k$ has a prime factor which 
occurs in ${\bf q}$ but not in ${\bf p}$, $k > j+1$ and 
$h(\widetilde{T}^n_k) \subset {\rm int}(T^n_m)$. 
We have therefore 
\[ w(T^n_m, T^n_0)=
w((T^n_m, h(\widetilde{T}^n_k))w(h(\widetilde{T}^n_k),
h(\widetilde{T}^n_j)) w(h(\widetilde{T}^n_j,T^n_0)).\] We have
obtained a contradiction because $q_k$ divides
$w(h(\widetilde{T}^n_k), h(\widetilde{T}^n_j))$ but not the left hand
side (which is non-zero also).  \end{proof}\end{proof}

\subsection{Open manifolds which are not w.g.s.c.}
\noindent In general the tower of obstructions we defined in the
previous sections is not trivial as is
shown below:

\begin{theorem} For uncountably many Whitehead-type manifolds $W^n$ 
of dimension $n\geq 3$ the manifolds  $W^n\times N^k$ are not 
$\infty$-compressible for any closed simply connected $k$-manifold
$N^k$.
\end{theorem}

\begin{proof}
It is sufficient to consider the case of the Whitehead-type manifolds since
the pair of groups appearing in the product exhaustions are the same
as this case.  

\vspace{0.2cm}
\noindent 
We start with the 3-dimensional case, and take for $W^3$ the classical
Whitehead manifold. Recall that $W^3$ is an
increasing union of solid tori $T_i$, with $T_i$ embedded in $T_{i+1}$
as a neighborhood of a Whitehead link. 

\vspace{0.2cm}
\noindent 
We shall first show that the pair $(T_{i+1},T_i)$ is not $2$-weakly
compressible, and hence not stably compressible. We then extend this
argument to show that any pair of the form $(T_{i+n},T_i)$ is not
$n+1$-weakly compressible, and hence not stably compressible. By
proposition~\ref{T:wgsc}, it follows that $W^3\times N^k$ is not
$\infty$-compressible, and hence not w.g.s.c.

\vspace{0.2cm}
\noindent 
Let $T$ and $T'$ be as usual and let $M^3=T'-{\rm int}(T)$, and fix a base
point $p\in\del T$. Then $C=G_0=\pi_1(M^3)$ in our usual notation.
Note that $\ker(\varphi)$ is normally generated by the meridian of $T$
and hence $\pi_1(M^3)/{\mathcal N}(\ker(\varphi),C)=\pi_1(T')={\mathbb
Z}$. Thus $G_1={\mathcal N}(\ker(\varphi),C)$ consists of the
homologically trivial elements in $\pi_1(M^3)$.

\vspace{0.2cm}
\noindent 
Consider now the cover $\widetilde{M^3}$ of $M^3$ with fundamental
group $G_1$. This is ${\mathbb R}^3$ with the neighborhood of an
infinite component link, say indexed by the integers, deleted. Further
each pair of adjacent components has linking number $1$. Pick a lift
$p'$ of the base point $p$, which we use for all the fundamental
groups we consider.

\vspace{0.2cm}
\noindent 
In this cover, $\psi(A)$ is the image of the bounding torus $T$ of the
component of this link containing $p'$, and $\psi(\ker(\varphi))$ is
generated by the meridian of this component. Thus, $G_1/{\mathcal
N}(\psi(\ker(\varphi), G_1)$ is the fundamental group of
$\widetilde{M^3}\cup_T D^2\times S^1$, i.e., of $\widetilde{M^3}$ with
a solid torus glued along $T$ to kill the meridian. But, because of
the linking, the longitude $\lambda\subset T$ is not trivial in this
group, i.e. $\lambda\notin G_2={\mathcal N}(\psi(\ker(\varphi),
G_1)$. Since $\lambda\in \psi(A)$, we see that the Whitehead link is
not $2$-compressible.

\vspace{0.2cm}
\noindent 
We shall now consider a pair $(T',T)$ in some refinement of the given
exhaustion. This is homeomorphic to a pair of the form $(T_n,T_1)$ for
some $n$. As before pick a base point $p\in\del T$.

\vspace{0.2cm}
\noindent 
Let $M_i=T_i-T_1$ and let $M=M_n$. Note that $M_i\subset M_{i+1}$. In
terms of earlier notation, $M_n=T'-T$ and $\pi_1(M)=C=G_0$. Further
$\ker(\varphi)$ is normally generated by the meridian of $T_1$.

\vspace{0.2cm}
\noindent 
We have a sequence of subgroups $G_k\subset G_0=\pi_1(M)$ and hence
covers $M^j$ of $M$ corresponding to these subgroups. Pick lifts $p^k$
of the base point $p$ to these covers. Then ${\mathcal
N}(\psi(\ker(\varphi), G_k)$ is generated by the meridian of the
component of the inverse image of $\del T_1$ containing $p^k$. As the
meridian is in $\ker(\varphi)$, and each $G_i$ is the normal subgroup
generated by $\ker(\varphi)$ in $G_{i-1}$, we see inductively that the
lift of the meridian is a closed curve in $M^k$ so that the previous
sentence makes sense.

\vspace{0.2cm}
\noindent 
Let $N^i$ be the result of gluing a solid torus or cylinder to $M^i$
along the component containing $p^k$ so that the meridian is
killed. Then by the above $G_k/{\mathcal N}(\psi(\ker(\varphi),
G_k)=\pi_1(N^i)$.

\vspace{0.2cm}
\noindent 
We shall prove by induction the following lemma.

\begin{lemma}
$T^{n-k}$ lifts to $N^k$, or equivalently, $M_{n-k}$ lifts to
$M^k$. Furthermore the longitude of the lift of $\del T_{n-k}$ is a
non-trivial element in $\pi_1(N^k)$.
\end{lemma}
\begin{proof}
The case when $k=1$ is the above special case. Suppose now that the
statement is true for $k$.

\vspace{0.2cm}
\noindent 
As the longitude of the lift of $\del T_{n-k}$ is a non-trivial
element in $\pi_1(N^k)$, in $M^{k+1}$ the inverse image of $M_{n-k}$
is a cylinder with a sequence of linked lifts of $M_{n-(k+1)}$
deleted. Thus, $M_{n-(k+1)}$ lifts to $M^{k+1}$, and its longitude is
linked with other lifts. It follows that the longitude of the lift of
$\del T_{n-(k+1)}$ is non-trivial in $\pi_1(N^{k+1})$.
\end{proof}

\noindent 
As a subgroup of $G_{n-1}$, $\psi(A)$ is the image of the lift of
$\del T_1$ containing the base point. As in the special case, as the
longitude of this torus is a non-trivial element of $\pi_1(N^{n-1})$,
it follows that $\psi(A)\not\subset G^n$. Thus $(T_n,T_1)$ is not
$n$-compressible.

\vspace{0.2cm}
\noindent 
This ends the proof of the claim for the Whitehead manifold. 
Observe however  that the same proof works for uncountably many similar manifolds --
namely we may embed $T_i$ in $T_{i+1}$ as a link similar to the
Whitehead link that winds around the solid torus several times.

\vspace{0.2cm}
\noindent 
We will use now a recurrence on the dimension and the results of the
previous section in order to settle the higher dimensional situation. 
Consider for simplicity $n=4$ and a  Whitehead-type manifold $W^4$ which is 
the ascending union of solid tori as in Wright's construction. We use
the notations from lemma \ref{wrconst} below.  Then the pair of tori $(T^4, Q^4)$ 
is constructed out of the two Whitehead links in one dimension less 
$(T_1^3, T_0^3)$ and $(D^2\times S^1, L^3)$. 

\vspace{0.2cm}
\noindent 
As above $\ker(\varphi)$ is normally generated by the meridian and 
$G_1$ consists of homologically trivial elements of
$\pi_1(T^4-{\rm int}(Q^4))$, by using Van Kampen and the fact that
$\pi_1(T^4)$ is abelian. 
The cover $\widetilde{M^4}$ of $M^4=T^4-{\rm int}(Q^4)$ with fundamental
group $G_1$ is ${\bf R}^4$ with a thick infinite link deleted. There
is an obvious ${\bf Z}^{2}$ action on the components of this link,
and so we can label the boundary tori as $T_{i,j}$, for integer $i,j$.

\vspace{0.2cm}
\noindent 
Let $\lambda$ be the longitude curve having the parameters 
$(n,k)$ on the torus $T_{0,0}$.  Then one can compute the linking numbers 
${\rm lk}(\lambda, T_{0,1})=k$ and ${\rm lk}(\lambda, T_{1,0})=n$. This follows
because both links used in the construction were the standard Whitehead
link. Variations which yield non-zero linking numbers are also
convenient for our purposes. Consequently for non-zero $n,k$ we 
obtained an element which is non-trivial in $G_2$ hence the pair of
solid tori is not $2$-compressible. 
A similar argument goes through the higher compressibility as
well. 
Using suitable variations in choosing the links and mixing the pairs 
of solid tori as in previous section yields uncountably many
examples as in the theorem. 
\end{proof}

\begin{remark}
It follows that $W^3\times D^k$ is not w.g.s.c. using the criterion
from \cite{Fun2,FT}. However the previous theorem is more precise
regarding the failure of g.s.c. for these product manifolds.
\end{remark}

\section{The proper homotopy invariance of the w.g.s.c.}

\subsection{Dehn exhaustibility}
\noindent We study in this section to what extent the w.g.s.c. is a proper
homotopy invariant.
\begin{definition}
A polyhedron $M$ is {\em (proper) homotopically dominated} by 
the polyhedron $X$ if there exists a map $f:M\to X$ such that 
the mapping cylinder $Z_f$ (properly) retracts on $M$. 
\end{definition}
\begin{remark}
A proper homotopy equivalence is the simplest example 
of a proper homotopically domination. 
\end{remark}

\noindent 
The main result of this section is:
\begin{theorem}\label{proper}
For $n\neq 4$ a non-compact  $n$-manifold is w.g.s.c. if and only if 
it is proper homotopically dominated by a w.g.s.c. polyhedron. 
\end{theorem}
\begin{remark}
It seems that the result does not hold, as stated,  for $n=4$
(see also the next section).
\end{remark}
\noindent {\em Proof of Theorem \ref{proper}}.
The main ingredient of the proof is the following 
notion, weaker than the w.g.s.c., introduced by Po\'enaru:
\begin{definition}
The simply-connected $n$-manifold  $W^n$ is {\em Dehn exhaustible} 
if, for any compact $K\subset W^n$ there exists some 
simply connected compact polyhedron $L$ and 
a commutative
diagram 
\[\begin{array}{ccc}
K &\; \; \;  \stackrel{f}{\to}& \; \; L \\
 & i\; \searrow & \; \;  \; \downarrow g \\
 & & \; \; W^n 
\end{array}
\]
\noindent 
where $i$ is the inclusion, $f$ is an embedding, $g$ is an 
immersion  and $f(K)\cap M_2(g)=\emptyset$.  
Here $M_2(g)$  is the set of double points, namely 
$M_2(g)=\{x\in L; \sharp g^{-1}(g(x))\geq 2\}\subset L$. 
If $n=3$ then one asks the map $g$ to be 
a {\em generic immersion}, which  means here that it has no triple
points.
\end{definition}
\noindent The first step is to establish:
\begin{proposition}
An open simply-connected manifold which is 
proper homotopically dominated by a w.g.s.c. polyhedron 
is Dehn exhaustible.
\end{proposition}
\begin{proof}
The proof given in \cite{Fun2} for  the 3-dimensional statement extends 
without any essential modification, and we skip the details. 
\end{proof} 
\begin{remark}
Po\'enaru proved a  Dehn-type lemma (see \cite{Po1}, p.333-339) which 
states that a Dehn exhaustible 3-manifold is w.g.s.c. This settles the 
dimension 3 case. 
\end{remark}
\begin{lemma}
If the open simply-connected $n$-manifold $W^n$ is Dehn exhaustible 
and $n\geq 5$ then it is w.g.s.c.
\end{lemma}
\begin{proof}
Consider a connected compact submanifold $K\subset W^n$. 
Assume that there exists a compact polyhedron $M^n$ with $\pi_1(M^n)=0$ and 
an immersion $F$ 
\[\begin{array}{ccc}
K &\; \; \;  \stackrel{f}{\to}& \; \; M \\
 & \; i\searrow & \; \;  \; \downarrow F \\
 & & \; \; W^n 
\end{array}
\]
\noindent such that $M_2(F)\cap K=\emptyset$.
\begin{lemma}\label{man}
One can suppose that $M^n$ is a manifold. 
\end{lemma}
\begin{proof}
The polyhedron $M^n$ is endowed  with an immersion
$F$ into the manifold
$W^n$. Among all abstract regular neighborhoods (i.e. thickenings) of   $M^n$
there is a $n$-dimensional one  $U(M^n,F)$, which is 
called the regular neighborhood determined by the immersion,  
such that the following conditions are fulfilled:
\begin{enumerate}
\item   $F:M^n\to W^n$ extends 
to an immersion $\widetilde{F}:U(M^n,F)\to W^n$.
\item The image of $\widetilde{F}(U(M^n, F))\subset W^n$ is the 
regular neighborhood of the polyhedron $F(M^n)$ in $W^n$. 
\end{enumerate}

\noindent  
The construction of the PL  regular neighborhood determined by an immersion 
of polyhedra is given in \cite{LiS}. The authors were building on the case 
of an  immersion of manifolds, considered previously  in \cite{HP}.  
Moreover, if one replaces  $M^n$ by the manifold $U(M^n, F)$ and 
$F$ by $\widetilde{F}$
we are in the conditions required by the Dehn-type lemma.
\end{proof}

\noindent Consider now  a handlebody decomposition of $M^n-f(K)$ 
and let $N_2^n$ be the 
union of $f(K)$ with the handles of index 1 and 2. Then $\pi_1(N_2^n)=0$. 
Let $\delta^2_j, \delta^1_j$ be the cores of these extra 1- and 2-handles. 
By using a small homotopy of $F$ one can replace
$F(\delta^2_j)\subset W^n$ by some embedded 2-disks $d_j^2\subset W^n$
with the same boundary. Also by general position these 2-disks can be chosen to
have disjoint interiors. Both assertions follow from the assumptions 
$n\geq 5$, and $M_2(F)\cap f(K)=\emptyset$. This implies that 
the restriction of the new map $F'$,  obtained by perturbing $F$, to 
$\delta^2_j$ (and $\delta^1_j$) is an embedding into $W^n-K$. 
Using the uniqueness of the regular neighborhood it follows that 
$F'$ can be chosen to be an embedding on $N_2^n$. In particular 
$K$ is engulfed in the 1-connected compact $F'(N_2^n)$.   
\end{proof}

\subsection{Dehn exhaustibility and end compressibility in dimension 4}
\begin{proposition}
An open 4-manifold is end compressible if and only if it  
is Dehn exhaustible.
\end{proposition}
\begin{proof}
We have to reconsider the proof of Proposition \ref{weak}. 
Everything works as above except that the disks 
$\delta^2_j$ cannot be anymore embedded, but only  (generically) immersed. 
They may have finitely many double points in their interior. 
Then the manifold $M^4$ obtained by adding 2-handles along the 
$\gamma_j$ has a generic immersion $F:M\to T'$, whose double  
points $M_2(F)$ are outside of $T$. This implies that 
$W^4$ is Dehn exhaustible. 

\vspace{0.2cm}
\noindent 
Conversely assume that $W^4$ is Dehn exhaustible. Let $K^4$ be a
compact submanifold of $W^4$ and $M^4$ be the immersible simply
connected polyhedron provided by the Dehn exhaustibility property.
Lemma \ref{man} allow us to assume that $M^4$ and $F(M^4)$ are
4-manifolds.  Consider now to the proof of the first claim from
Theorem \ref{proper}. It is sufficient to consider the case when
$M^4=M_2^4$ i.e. $M^4$ is obtained from $K^4$ by adding 1- and
2-handles. If $\Gamma\subset \pi_1(\partial M_1^4)$ is the normal
subgroup generated by the attaching curves of the 2-handles of $M^4$
then the same argument yields:
\[ \Gamma \subset \ker(\pi_1(\partial M_1^4)\to \pi_1(M^4-{\rm int}(K^4))).\]
Since $F$ is a generic immersion we can suppose that $F$ is an 
embedding of the cores of the 1-handles and so $F|_{M_1^4}$ is an 
embedding. 

\vspace{0.2cm}
\noindent 
We have $F(M^4-{\rm int}(K^4))\subset F(M^4)-{\rm int}(K^4)$ because the double points 
of $F$ are outside $K^4$. Now the homomorphism induced by $F$ on the
left side of the diagram 
\[\begin{array}{cccc}
\pi_1(\partial M_1^4) & \to & \pi_1(M^4-{\rm int}(K^4)) \\
F \downarrow         &     & \downarrow F          \\
\pi_1(\partial F(M_1^4))& \to & \pi_1(F(M^4)-{\rm int}(K^4))
\end{array}\]
\noindent 
is an isomorphism and we derive that 
\[F(\Gamma) \subset \ker(\pi_1(\partial F(M_1^4))\to \pi_1(F(M^4)-{\rm int}(K^4))).\]

\noindent 
Meantime $F(\Gamma)$ surjects onto $\pi_1(F(M_1^4))$ under the map
$\pi_1(\partial F(M_1^4))\to \pi_1(F(M_1^4))$.  But $F(M_1^4)$ is
homeomorphic to $M_1^4$, hence it is obtained from $K^4$ by adding
1-handles.  This shows that the pair $(F(M^4), K^4)$ is stably
compressible, from which one obtains the end compressibility of $W^4$
as in the proof of Theorem \ref{wgsc}.  \end{proof}
\begin{remark}
If the open 4-manifold $W^4$ is Dehn exhaustible then $W^4\times
[0,1]$ is also Dehn exhaustible hence w.g.s.c. Therefore an example of
an open 4-manifold $W^4$ which is end compressible but which is not
w.g.s.c.  will show that the result of Theorem \ref{proper} cannot be
extended to dimension 4, as stated.  Such examples are very likely to
exist, as the Dehn lemma is known to fail in dimension 4 (by
S.~Akbulut's examples).
\end{remark}

\subsection{Proper-homotopy invariance of the end compressibility}
\noindent It would be interesting to have a soft version of the theorem
\ref{proper} for the end compressibility situation. Notice that the
definition of the end compressibility extends word by word to
non-compact polyhedra. One uses instead of the boundary of manifolds
the frontier of a polyhedron.
\begin{remark}
One expects that the following be true. 
If a polyhedron $M$ is proper homotopically dominated by an end compressible 
polyhedron $X$ then $M$ is also end compressible. 
The only ingredient lacking for the complete proof is the 
analogue of the remark \ref{indep} for polyhedra: if one exhaustion is
stably-compressible then all exhaustions have stably-compressible 
exhaustions. 
\end{remark}
\noindent Along the same lines  we have: 
\begin{proposition}
If there is a degree one map $f:X^n\to M^n$ between one ended manifolds of
the same dimension, then if $X^n$ is end-compressible so is $M^n$.
\end{proposition}
\begin{proof}
 
\noindent We use the fact that  degree-one maps are surjective on fundamental
group. Given an exhaustion $\{L_j\}$ of $M^n$, pull it back to
$\{f^{-1}(L_j)\}$ of $X^n$. Notice that $\del f^{-1}(L_j)=f^{-1}(\del L_j)$
where $\del$ stands for the frontier. 

\vspace{0.2cm}
\noindent 
One needs then the following approximation
by manifolds result. Given two $n$-complexes $K_1\subset
{\rm int}(K_2)\subset {\bf R}^n$ there exist regular neighborhoods 
$K_j^{\varepsilon}\subset{\bf R}^n$ such that 
$K_j\subset {\rm int}(K_j^{\varepsilon})$, 
$K_j$ is homotopy equivalent to $K_j^{\varepsilon}$, 
$K_2-{\rm int}(K_1)$ is homotopy equivalent to
$K_2^{\varepsilon}-K_1^{\varepsilon}$, and moreover 
$\partial K_1^{\varepsilon}$ is homotopy equivalent to $\partial
K_1$. This uses essentially the fact that $K_j$ are of codimension
zero in ${\bf R}^n$. 

\vspace{0.2cm}
\noindent  
Now, the hypothesis applied to  the approximating exhaustion
consisting of submanifolds implies the existence of a stably-compressible refinement
of $\{f^{-1}(L_j)\}$. Since degree-one maps are surjective 
on fundamental groups the lemma below permits to descend to $M^n$.  

\begin{lemma}
Suppose that  the triple $(A,B,C)$ of groups surjects onto $(A',B',C')$ 
i.e. we have three surjections with diagrams commuting. Then if
$(A,B,C)$ is strongly (or stably) compressible then so is
$(A',B',C')$. 
\end{lemma}
\noindent The proof is straightforward.

\vspace{0.2cm}
\noindent
The only subtlety above is to make sure the inverse image of
boundary components is connected (else we can connect them up in the
one-ended case).
\end{proof}
        
\begin{remark}
In the many-ended case, we need to say that we have a degree-one map
between each pair of ends (not 2 ends mapping to one, with one of them
having degree 2 and the other -1). This holds in particular for a
proper map that has degree one and is injective on ends.
\end{remark}

\section{G.s.c. for 4-manifolds}
\subsection{W.g.s.c. versus g.s.c.}
\begin{definition}
A {\em geometric Po\'enaru-Mazur-type manifold} $M^4$ is a 
compact simply connected  $4$-manifold satisfying the following conditions:
\begin{enumerate}
\item $H_2(M^4)=0$.
\item the boundary $\partial M^4$ is connected and 
$\pi_1$-dominates a virtually geometric 3-manifold group, i.e.
there exists a surjective homomorphism
\[ \pi_1(\partial M^4)\to \pi_1(N^3), \]
onto the (non-trivial) fundamental group of a virtually 
geometric 3-manifold $N^3$.
\end{enumerate}
\end{definition}
\begin{proposition}\label{notgsc}
The interior ${\rm int}(M^4)$ of a geometric Po\'enaru-Mazur-type manifold
$M^4$ does not have a proper handlebody decomposition without
1-handles with the boundary of a cofinal subset of the
intermediate manifolds obtained on a finite number of handle additions
being homology spheres.
\end{proposition}
\subsection{Casson's proof of Proposition \ref{notgsc}}
The main ingredient is the following proposition
extending an unpublished result of A.~Casson:
\begin{proposition}\label{T:Casson}
Consider  the 4-dimensional (compact) cobordism $(W^4, M^3, N^3)$ such that 
$(W^4, M^3)$ is 1-connected. Assume moreover that the following conditions 
are satisfied: 
\begin{enumerate}
\item $H_2(W^4, M^3; {\bf Q}) = 0$, both $M^3$ and $N^3$ are 
connected. 
\item $\pi_1(N^3)$ is a group which $\pi_1$-dominates a virtually 
geometric non-trivial 3-manifold group. Let $K$ be the kernel 
of this epimorphism.
\item $b_1(W^4)\leq b_1(N^3)$, where $b_1$ denotes the first Betti number.
\item The map $\pi_1(N^3)\to \pi_1(W^4)$ induced by the 
inclusion $N^3\hookrightarrow W^4$ has kernel strictly bigger than 
the subgroup $K$. In particular this is true if this map is trivial. 
\end{enumerate} 
Then any handlebody decomposition of $W^4$ from $M^3$ has 1-handles
i.e. the pair $(W^4, M^3)$ is not g.s.c.
\end{proposition}
\begin{remark}
A necessary condition for the g.s.c. of 
$(W^4, M^3)$ is that the map $\pi_1(M^3)\to \pi_1(W^4)$, induced by the 
inclusion $M^3\hookrightarrow W^4$, be onto. In fact adding 2-handles 
amounts to introducing new relations to the fundamental group of the 
boundary, whereas the latter is not affected by higher dimensional handle
additions.
\end{remark}

\noindent 
Casson's result was based on partial positive solutions to the
Kervaire Conjecture \ref{kerv}. 
One proves that certain $4$-manifolds $(N,\del N)$ have no handle
decompositions without $1$-handles by showing that if they did, then
$\pi_1(\del N)$ violates the Kervaire conjecture.
Casson's argument works to the extent that the Kervaire conjecture is
known to be true. Casson originally applied it using a theorem of
M.~Gerstenhaber and O.S.~Rothaus (\cite {GR}), which said that the Kervaire
conjecture holds for subgroups of a compact Lie group. Subsequently,
O.S.~Rothaus (\cite{Rot}) showed that the conjecture in fact holds for
residually finite groups. Since residual finiteness for all
$3$-manifold groups is implied by the geometrization conjecture,
Casson's argument works in particular for all manifolds satisfying the
geometrization conjecture. A simple argument (Remark~\ref{T:alg}
below) extends the class of groups for which the Kervaire conjecture
is known further.

\begin{proposition} \label{T:alg} If some non-trivial quotient $Q$ of a
group $G$ satisfies the Kervaire conjecture, then so does $G$. In 
particular if a finitely generated group $G$ has a proper
finite-index subgroup, then G satisfies the Kervaire conjecture
(since finite groups satisfy the Kervaire conjecture by 
\cite{GR}).
\end{proposition}

\begin{proof} Let $\phi:G\onto Q$ be the quotient map. Assume that $Q$
satisfies the Kervaire conjecture. Suppose that G violates the
Kervaire conjecture. Then we have generators $\alpha_1,\dots,\alpha_n$
and relations such that
$\frac{G*\langle\alpha_1,\dots,\alpha_n\rangle}{
\langle\langle r_1,\dots,r_n\rangle\rangle}$ is the trivial
group. Let
$\phi^*:G*\langle\alpha_1,\dots,\alpha_n\rangle\to Q*\langle
\overline{\alpha_1},\dots,\overline{\alpha_n}\rangle$ be the map extending $\phi$ by mapping
$\alpha_i$ to $\overline{\alpha_i}$. This is clearly a surjection, and
induces a surjective map
$\overline{\phi}:\frac{G*\langle\alpha_1,\dots,\alpha_n\rangle}{
\langle\langle r_1,\dots,r_n\rangle\rangle}\onto\frac{Q*\langle 
\overline{\alpha_1},\dots,\overline{\alpha_n}\rangle}{\langle\langle
\phi^*(r_1),\dots,\phi^*(r_n)\rangle\rangle}$. But
since the domain of the surjection $\overline{\phi}$ is trivial, so is the
codomain. But this means that
$\frac{Q*\langle
\overline{\alpha_1},\dots,\overline{\alpha_n}\rangle}
{\langle\langle\phi^*(r_1),\dots,\phi^*(r_n)\rangle\rangle}$
is trivial, and so $Q$ violates the Kervaire conjecture, a
contradiction.
\end{proof}

\begin{proof}[Proof of Proposition~\ref{T:Casson}] 
Suppose that 
\[ W^4=M^3\times[0,1]\sharp_k 2{\mbox{\rm -handles} } \sharp_r 3{\mbox 
{\rm -handles}},\]
(with some 0-handle or 4-handle added if one boundary component is 
empty). 
It is well-known (see \cite{RS}) that the homology groups
$H_*(W^4, M^3)$ are the same as those of a differential complex 
$C_*$, whose component $C_j$ is the free module generated by the 
$j$-handles. Therefore this complex has the form:
\[ 0\to {\bf Z}^r\to {\bf Z}^k\to 0.\]
Thus  $H_2(W^4, M^3;{\bf Q})=0$ implies  that $k\leq r$ holds. 

\vspace{0.2cm}
\noindent 
Consider now the handlebody decomposition is turned up-side down:
\[W^4=N^3\times[0,1]\sharp_r 1{\mbox{\rm -handles }} \sharp_k 2{\mbox{\rm -handles}},\] 
(plus possibly one 0-handle or 4-handle if the respective boundary 
component is empty).

\vspace{0.2cm}
\noindent 
By the van Kampen theorem it follows that 
$\pi_1(W^4)$ is obtained from $\pi_1(N^3)=\pi_1(N^3\times[0,1])$ 
by adding one generator for each 1-handle and one relation for each 
2-handle. 
Therefore 
\[ \pi_1(W^4) = \pi_1(N^3) * {\bf F}(r)/W(k), \]
where ${\bf F}(r)$ is the free group on $r$ generators $x_1,...,x_r$ and 
$W(k)$ is a normal subgroup of the free product generated also 
by $k$ words $Y_1,...,Y_k$.

\vspace{0.2cm}
\noindent 
Consider a virtually geometric 3-manifold $L^3$ such that 
$\pi_1(M^3)\to \pi_1(L^3)$ is surjective. 
If $L^3$ is a geometric 3-manifold  then its fundamental group is 
residually finite (see e.g. \cite{Th}, Thm.3.3, p.364). 
Let $d_{ij}$ be the degree of the letter $x_j$ in the word 
representing $Y_i$. 
The result of Rothaus (\cite{Rot}, Thm. 18, p.611) states that 
for any locally residually finite group $G$ and choice of words 
$Y_i$ such that ${\bf d}=(d_{ij})_{i,j}$ is of  (maximal) rank $k$,
 the natural morphism $G \to G * {\bf F}(r)/W(k)$ is an injection. 
We have therefore a commutative diagram
\[\begin{array}{ccc}
\pi_1(M^3) & \to & \pi_1(W^4) \\
\downarrow &     & \downarrow \\
 \pi_1(L^3)& \hookrightarrow   &\pi_1(L^3)* {\bf F}(r)/W(k)
\end{array}\]
whose vertical arrows are surjections. 
The kernel of the map induced by inclusion, 
$\pi_1(M^3) \to  \pi_1(W^4)$ is contained in $K$. 
This contradicts our hypothesis.

\vspace{0.2cm}
\noindent 
On the other hand if the rank of ${\bf d}$ is not maximal then by considering 
the abelianisations one derives 
$H_1(G * {\bf F}(r)/W(k))\subset  H_1(G) \oplus {\bf Z}$, 
hence $b_1(W^4)\geq b_1(N^3)+1$, which is also false.
\end{proof}

\begin{corollary}
Consider a 4-manifold $W^4$ which is compact connected simply-connected 
with non-simply connected boundary $\partial M$. If the boundary is 
(virtually) geometric 
and $H_2(W^4)=0$ then $W^4$ is not g.s.c. 
\end{corollary}

\begin{proof}[Proof of Proposition \ref{notgsc}]
Assume now that ${\rm int}(M^4)$ admitted a proper handlebody decomposition 
without 1-handles. One identifies ${\rm int}(M^4)$ with 
$M^4\cup_{\partial M\cong \partial M\times\{0\}} \partial M\times [0,1)$.
We can truncate the handle decomposition at a finite stage 
in order to obtain a manifold $Q^4$ such that 
$\partial Q^4 \subset \partial M^4 \times (0,1)$, because the decomposition
is proper. We can suppose that $\partial Q^4$ is connected since 
${\rm int}(M^4)$ has one end. Then $Q^4$ is g.s.c. hence $\pi_1(Q^4)=0$. 

\vspace{0.2cm}
\noindent 
By hypothesis, we can choose   $\partial Q^4$ to be a homology 
sphere.  Then $\partial Q^4$ separates  the cylinder $\partial M^4\times[0,1]$ 
into two manifolds with boundary which, by Mayer-Vietoris, have 
the homology of $S^3$. This implies that
$H_2(Q^4)=0$ (again by Mayer-Vietoris). 

\vspace{0.2cm}
\noindent 
Let us consider now the map $f:\partial Q^4\hookrightarrow \partial M^4\times[0,1]
\to \partial M^4$, the composition of the inclusion with the obvious 
projection. 

\begin{lemma}
The map $f$ has degree one hence induces a surjection on the fundamental 
groups.
\end{lemma}
\begin{proof}
The 3-manifold $\partial Q^4$  separates the two components 
of the boundary. In particular the generic arc joining 
$\partial M^4\times\{0\}$ to 
$\partial M^4\times\{1\}$ intersects transversally $\partial Q^4$ 
in a number of points, which counted with the sign sum up to 1 (or -1).
If properly interpreted this is the same as claiming the degree of $f$ is one. 

\vspace{0.2cm}
\noindent 
It is well-known that a degree one map between orientable 3-manifolds 
induce a surjective map on the fundamental group (more generally, 
the image of the homomorphism induced by a  degree $d$ is a 
subgroup whose index is bounded by $d$). \end{proof}

\vspace{0.2cm}
\noindent 
This shows that $\pi_1(\partial Q^4)\to \pi_1(\partial M^4)$ is 
surjective. On the other hand $\pi_1(\partial M^4)$ surjects onto 
a non-trivial residually finite group. 
Since $\pi_1(\partial Q^4)\to \pi_1(Q^4)=1$ is the trivial map, the 
argument we used previously (from Rothaus' theorem) 
gives us a contradiction. This settles our claim.
\end{proof}

\section{Handle decompositions without 1-handles in dimension 4}
\subsection{Open tame 4-manifolds}
\begin{definition}
An exhaustion of a 4-manifold is {\em g.s.c.} 
if it corresponds to a proper sequence of
handle additions with no 1-handles.  Alternatively one has a proper
Morse function, which we will refer to as {\em time}, with words like
{\em past} and {\em future} having obvious meanings, with no critical
points of index one. The inverse images of regular points are
$3$-manifolds, which we refer to as the {\em manifold at that time}.
\end{definition}

\noindent We assume henceforth that we have a g.s.c. handle decomposition of the
interior ${\rm int}(W^4)$ of $(W^4,\del W^4)$, a compact four manifold with
boundary a homology $3$-sphere and $H_2(W^4)=0$.

\vspace{0.2cm}
\noindent 
Now let $(K_i^4,\del K_i^4), i\in \mathbb N$ denote the $4$-manifolds
obtained by successively attaching handles to the zero handle
$(B^4,S^3)$, that is if $t:(W^4,\del W^4)\to \mathbb R$ is the Morse
function time, then $(K_i^4,\del K_i^4) = t^{-1}((\infty,a_i])$, with
$a_i$ being points lying between pairs of critical values of the Morse
function.

\begin{lemma}\label{T:surgery} $\del K_{i+1}^4$ is obtained from $\del K_i^4$ by one of the following:
\begin{itemize}
\item A 0-frame surgery about a homologically trivial knot in $\del
K_i^4$.
\item Cutting along a non-separating $2$-sphere in $\del K_i^4$ and
capping off the result by attaching a $3$-ball.
\end{itemize}
These correspond respectively to attaching $2$-handles and $3$-handles
to $(K_i^4,\del K_i^4)$.
\end{lemma}
\begin{proof} Since attaching $2$-handles and $3$-handles 
correspond to surgery and cutting along $2$-spheres respectively, we
merely have to show that the surgery is 0-frame about a homologically
trivial curve and the spheres along which one cuts are non-separating.

\vspace{0.2cm}
\noindent 
First note that the absence of  $1$-handles implies
$H_1(K_i^4)=0=\pi_1(K_i^4)$, for all $i$. Further, each $\del K_i^4$ is
connected because  ${\rm int}(W^4)$ has one end. 
Thus the $2$-spheres along which any $\del K_i^4$ is
split have to be non-separating.

\vspace{0.2cm}
\noindent Using Mayer-Vietoris, the fact that $H_2(W^4)=0$, and the long exact
sequence in homology we derive that $H_2(K_i^4)=H_2(\del K_i^4)$. Also 
adding a $3$-handle decreases the rank
of $H_2(\del K_i^4)$ by one hence every surgery 
increases the rank of $H_2(\del K_i^4)$ by one unit. But this means that the surgery must be a zero-frame surgery
about a homologically trivial curve.
\end{proof}

\noindent 
For $i$ large
enough, $\del K_i^4$ lies in a collar $\del W^4\times [0,\infty)$ 
hence we have a map $f_i:\del
K_i^4\to\del W^4$ which is the composition of the inclusion with the
projection. By the Lemma 5.1 the maps $f_i$ are of degree one and 
induce surjections $\phi_i:\pi_1(\del K_i^4)\to \pi_1(\del W^4)$. Here and
henceforth we always assume that the index $i$ is large enough so that
$f_i$ is defined.

\begin{lemma}\label{T:degreeone} 
The homotopy class of a curve along which surgery
is performed is in the kernel of $\phi_i:\pi_1(\del K_i^4)\to \pi_1(\del W^4)$.
\end{lemma}
\begin{proof} If a surgery is performed along a curve $\gamma$, this means
that a 2-handle is attached along the curve in the
$4$-manifold $W^4$. Hence $\gamma$ bounds a disk in $\del
W^4\times[0,\infty)$, which projects to a disk bounded by $f_i(\gamma)$
in $\del W^4$.
\end{proof}

\begin{remark}\label{newmap} The maps $\phi_i$ and $\phi_{i+1}$ are
related in a natural way. To define the map $\phi_{i+1}$, take a
generic curve $\gamma$ representing any given element of $\pi_1(\del
K_{i+1}^4)$. If $\del K_{i+1}^4$ is obtained from $\del K_i^4$ by splitting
along a sphere, then $\gamma$ is a curve in $\del K_i^4$, and so we can
simply take its image. On the other hand, if a surgery was performed,
then we may assume that $\gamma$ lies off the solid torus that has
been attached, and hence lies in $\del K_i^4$, so we can take its image
as before. This map is well-defined by lemma~\ref{T:degreeone}.
\end{remark}

\begin{definition}
 A curve $\gamma'\subset \partial K_i^4$ is a {\em descendant} of the 
surgery curve $\gamma \subset \partial K_i^4$ if it is homotopic to 
it in $\partial K_i^4$ (though not in general homotopic to $\gamma$ after
the surgery). 
A curve $\gamma\subset \del K_i^4$ is said to {\em persist till
$\del K_{i+n}^4$} if some descendant of $\gamma$ persists, i.e., we can
homotope $\gamma$ in $\del K_i^4$ so that it is disjoint from all the future
$2$-spheres on which 3-handles are attached while passing from 
$\del K_i^4=M_i^3$ to $\del K_{i+n}^4=M_{i+n}^3$.
\end{definition}

\begin{definition} A curve $\gamma\subset \del K_i^4$ is said to {\em die by
$\del K_{i+n}^4$} if it is homotopically trivial in the $4$-manifold
obtained by attaching $2$-handles to $K_i^4$ along the curves in $\del
K_i^4$ where surgeries are performed in the process of passing to $\del
K_{i+n}^4$, or equivalently, $\gamma$ is trivial in the group obtained
by adding relations to $\pi_1(M_i^3)$ corresponding to curves along
which the surgery is performed.
\end{definition}

\noindent 
We prove now a key property of the sequence $\del K_i^4$.
\begin{lemma}\label{T:persistsdies} For each $i$, there is a uniform 
$n=n(i)$ such
that any curve $\gamma\subset \del K_i^4$, 
$\gamma\in \ker \phi_i$
that persists till $\del K_{i+n}^4$ dies by $\del K_{i+n}^4$.
\end{lemma}
\begin{proof} We can
find $x\in [0,\infty)$ so that $\del W^4\times \{x\}$ is entirely
after $\del K_i^4$, and $n_1$ so that $\del W^4\times \{x\}\subset
K_{i+n_1^4}$, because the handlebody decomposition is proper.  We then
define $n$ by repeating this process once, i.e. $\del W^4\times
\{x_1+\varepsilon\}\subset K_{i+n}^4$, for some $x_1+\varepsilon >x_1 >
x$ for which $\del W^4\times \{x_1\}$ is entirely after $\del
K_{i+n_1}^4$.  Consider $\gamma\in \ker\phi_i$ which persists till $\del
K_{i+n_1}^4$. This means that there is an annulus properly embedded in
$K_{i+n}^4-{\rm int}(K_i^4)$, whose boundary curves are $\gamma$ and $\tilde
\gamma \in \del K_{i+n}^4\subset \del W^4\times [x_1,
x_1+\varepsilon)$. Since $\gamma \in \ker\phi_i$ it bounds a disc in
$\del W^4\times [x_1,\infty)$. This disc together with the above
annulus ensure that $\gamma$ dies by $\del K_{i+n}^4$, as they bound
together a disc entirely in $K_{i+n}^4-{\rm int}(K_i^4)$, and $3$-handles do not
affect the fundamental group.  \end{proof}

\subsection{The structure theorem}\label{S:reorder}

\noindent 
Suppose henceforth that we have a
sequence of  connected $3$-manifolds $M_i^3\subset \del W^4\times [0,\infty)$ 
and associated maps onto $f_i:M_i^3\to \del W^4$
that satisfies the properties of $\del K_i$ stated above. Specifically
one asks that: 
\begin{itemize}
\item The maps $f_i:M_i^3\to \del W^4$ are of degree one, hence inducing 
surjection $\phi_i:\pi_1(M_i^3)\to \pi_1(\del W^4)$. 
\item $M_{i+1}$ is obtained from $M_i^3$ either by a 0-frame surgery 
along a homologically trivial knot in $M_i^3$, or else by cutting 
along a non-separating 2-sphere in $M_i^3$. 
\item The surgery curves in $M_i^3$ belong to $\ker\phi_i$.
\item The maps $\phi_i$ and $\phi_{i+1}$ are related as in 
Remark~\ref{newmap}. 
\item For any $i$ there exists some $n=n(i)$ such that 
any curve in $M_i^3$ which persists till $M_{i+n}^3$ dies by $M_{i+n}^3$.
\end{itemize}

\noindent 
We show in this section that, after possibly changing the order of
attaching handles, any handle decomposition without $1$-handles is of
a particular form.

\vspace{0.2cm}
\noindent 
We first describe a procedure for attempting to construct a handle
decomposition for ${\rm int}(W^4)$ starting with a partial handle
decomposition, with boundary $M_i^3$. In general, $M_i^3$ has
non-trivial homology. It follows readily from the proof of
lemma~\ref{T:surgery} that $H_1(M_i^3)$ is a torsion free abelian
group. The only way we can remove homology is by splitting along
spheres. To this end, we take a collection of surfaces representing
the homology, perform surgeries along curves in these surfaces so that
they compress down to spheres, and then split along these spheres.  By
doing the surgeries, we have created new homology, and hence have to
take new surfaces representing this homology and continue this
procedure. In addition to this, we may need to perform other surgeries
to get rid of the {\em homologically trivial portion} of the kernel of
$\phi_i:\pi_1(M_i^3) \to \pi_1(\del W^4)$.

\vspace{0.2cm}
\noindent 
The above construction may meet obstructions, since the surgeries have
to be performed about curves that are homologically trivial as well as
lie in the kernel of $\phi_i$, hence it may not be always possible to
perform enough of them to compress the surfaces to spheres. The
construction terminates at some finite stage if at that stage all the
homology is represented by spheres and no surgery off these surfaces
is necessary.

\vspace{0.2cm}
\noindent 
\begin{theorem}\label{T:reorder} After possibly changing the order of
attaching handles, any handle decomposition without $1$-handles may be
described as follows. We have a collection of surfaces
$F_j^2(i)$, with disjoint
simple closed curves $l_{j,k}\subset F_j^2(i)$ and
a generic  immersion $\psi_i: \cup_j F_j^2(i) \to M_i^3$ such that:
\begin{itemize}
\item The surfaces represent  the homology of $M_i^3$, i.e. 
$\psi_i$ induces a surjection $\psi_i:H_2(\cup_j F_j^2(i))\to H_2(M_i^3)$.
\item  The immersion $\psi_i$ has only ordinary double points 
and the restriction to each individual surface $F_j(i)$ is an
embedding. The double curves of $\psi_i$ are among the curves $l_{j,k}$. 
Their images $\psi_i(l_{j,k})$ are called seams. 
\item When compressed along the seams (i.e. by adding 2-handles 
along them) the surfaces $\psi_i(F_j^2(i))$  become unions of spheres. 
\item The seams are homologically trivial curves in $M_i^3$ and lie in
the kernel of $\phi_i$.
\item The pull backs (see the definition 
below)  of the surfaces  $\psi_m(F_k^2(m))\subset M_m^3$ for 
$m > i$, which are surfaces with boundary in $M_i^3$, can only intersect the
$F_j(i)$'s either transversely at the seams or by having some
boundary components along the seams.
\end{itemize}
We attach $2$-handles along all the seams of $M_i^3$, and possibly also
along some curves that are completely off the surfaces
$F_j^2(i)$ in $M_i^3$ and
have no intersection with any future surface $F_k(m)$, $m > i$. 
We then attach
$3$-handles along the 2-spheres obtained by compressing the
surfaces $\psi_i(F_j^2(i))$. 
Iterating this procedure gives us the handle decomposition.
\end{theorem}

\noindent 
We will  see that once we construct the surfaces, all of the
properties follow automatically. 

\vspace{0.2cm}
\noindent 
Let $F^2\subset M_{i+n}^3$ be an embedded surface. 
We let $M_i^3=\del K_i^4$, where $K_i^4$ is the bounded 
component in ${\rm int}(W^4)$. Then 
$K_{i+n}^4-{\rm int}(K_i^4)= M_i^3\times [0,\varepsilon]\cup h_j^2\cup h_k^3$, 
where $h_j^m$ are the  attached $m$-handles.
 
\begin{lemma}
There exists an isotopy of $K_{i+n}^4-{\rm int}(K_i^4)$ such that 
$F^2\subset M_i^3\times \{\varepsilon\} \cup h_j^2$ and 
$F^2 \cap h_j^2 = \cup_k \delta^2_{j,k}$, where 
$\delta^2_{j,k}$ are disjoint 2-disks properly embedded in the 
pair $(h_j^2,\partial_a h_j^2)$ (here $\partial_a h_j^2$ denotes 
the attachment zone of the handle, which is a solid torus), 
which are parallel to the core
of the handle. Moreover $\partial \delta^2_{j,k}\subset \partial(\partial_a
h_j^2))$ are concentric circles on the torus, parallel to the 
0-framing of the attaching circle. 
\end{lemma}
\begin{proof}
It follows from a transversality argument that the image of $F$ 
intersects only the 2-handles, along 2-disks. 
Further it is sufficient to see that the 
circles  $\partial \delta^2_{j,k}$ are homotopic to the 0-framing 
since in $K_{i+n}^4-{\rm int}(K_i^4)$ homotopy implies isotopy for 
circles. If one circle is null-homotopic then it can be removed by
means of an ambient isotopy. If a circle turns
$p$-times around the longitude, then it cannot bound a disk in $h_j^2$
unless $p=1$. 
\end{proof}
\begin{definition}
Consider a parallel copy in $M_i^3=M_i^3\times\{0\}$ of the surface
with boundary $F'^2= F^2- \cup_{j,k}\delta^2_{j,k}\subset
M_i^3\times\{\varepsilon\}- \cup_j \partial_a h_j^2$, and use
standardly embedded annuli in the torus $\partial_a h_j^2$, which join
the parallel circles to the central knot in order to get a surface
with boundary on the surgery loci. We calls this a {\em pull back} of
the surface $F^2\subset M_{i+n}^3$.
\end{definition}

\begin{lemma}\label{T:homonto} Let $\alpha:\pi_1(M_i^3)\to H_1(M_i^3)$ be
the Hurewicz map. Then $\phi_i(\ker(\alpha))=\pi_1(\del W^4)$, i.e. 
the pair $(\alpha, \phi_i)$
is strongly compressible.
\end{lemma}
\begin{proof} Consider the diagram 
$\begin{array}{ccc}
\Gamma & \stackrel{\phi}{\longrightarrow} & G \\
\pi \downarrow &  &  \downarrow \\
\Gamma_{ab} & \stackrel{\phi_{ab}}{\longrightarrow} & G_{ab}
\end{array}$
where $\phi$ is surjective, and the subscript $ab$ means
abelianisation. Then it is automatically that 
$\pi(\ker\phi)=\ker\phi_{ab}$. 
Since $H_1(\del W^4)=0$, and the strong compressibility is 
symmetric, the result follows.
\end{proof}

\noindent 
Now, let $n=n(i)$ be as in the conclusion of
lemma~\ref{T:persistsdies}. We consider a maximal set of disjoint 
non-parallel essential 2-spheres (which is uniquely defined up to isotopy) 
and pull back these spheres up to time $i$ to get a collection of planar
surfaces, whose union is a 2-dimensional polyhedron $\Sigma_i\subset M_i^3$. 
\begin{lemma}
If  $\iota$ denotes the map induced by the 
inclusion $\pi_1(M_i^3-\Sigma_i)\to \pi_1(M_i^3)$ then  the restriction
\[ \phi_i:\iota(\pi_1(M_i^3-\Sigma_i))\cap \ker(\alpha)
\to \pi_1(\del W^4)\] 
is surjective. 
\end{lemma}
\begin{proof} 
The pull-backs in $M_i^3$ of  spheres $S^2_m\subset M_{i+j}^3$ 
are planar surfaces with boundary components being the loci of future
surgeries. Further, after compressing the spheres $S^2_m$  of $M_{i+j}$ 
(hence arriving into $M_{i+j+k}^3$) we have 
a surjection  $\phi_{i+j+k}$, thus the map 
 $\pi_1(M_{i+j}^3-\cup S^2_m)\to \pi_1(\del W^4)$ is also surjective. 
This means that there exist
curves in the complement of the planar surfaces in $M_i^3$  mapping to every
element of $\pi_1(\del W^4)$. 
Moreover, by the above lemma, we have such curves that are
homologically trivial in $M_{i+j}^3$, and hence in $M_i^3$ as all
surgery curves are null-homologous. 
\end{proof}
\begin{lemma}\label{T:homcutup} $i_*:H_1(M_i^3 - \Sigma_i)\to
H_1(M_i^3)$ is the zero map.
\end{lemma}
\begin{proof} 
If not then there exists a curve $\gamma\subset M_i^3 - \Sigma_i$ that
represents a non-trivial element of $H_1(M_i^3)$.  Modifying by a
homologically trivial element if necessary, we may assume that
$\gamma\in ker(\phi_i)$. By the previous lemma $\gamma$ persists. The
group $\pi_1(K_{i+n}^4-{\rm int}(K_i^4))$ is the quotient of $\pi_1(M_i^3)$ by
the relations generated by the surgery curves, which are homologically
trivial. In particular $H_1(K_{i+n}^4-{\rm int}(K_i^4))=H_1(M_i^3)$. Then the
class of $\gamma\in H_1(K_{i+n}^4-{\rm int}(K_i^4))$ is non-zero since its image
in $H_1(M_i^3)$ is non-zero by hypothesis.  This gives the required
contradiction.  \end{proof}

\noindent 
We are now in a position to prove the structure theorem. The
images of the immersion $\psi_i$ is 
obtained from the polyhedron $\Sigma_i$ by {\em stitching together}
several  planar surfaces  along the boundary knots. 
These knots will be the seams of the
surfaces. It is clear by construction that we have all the desired
properties as soon as we show that there are enough planar surfaces to
be stitched together to represent all the homology.
 
\vspace{0.2cm}
\noindent        
To see this, we consider the reduced homology exact-sequence of the
pair $(M_i^3,M_i^3 - \Sigma_i)$, and use the fact that
$M_i^3 - \Sigma$ is connected, since $M_{i+n}^3$ is, as well as
lemma~\ref{T:homcutup}. Thus, we have the exact sequence
$$\dots\to H_1(M_i^3 - \Sigma_i)\to H_1(M_i^3)\to H_1(M_i^3,
M_i^3 - \Sigma_i)\to\tilde H_0( M_i^3 - \Sigma_i)$$ which gives
the exact sequence
$$0\to H_1(M_i^3) \to H_1(M_i^3, M_i^3 - \Sigma_i)\to 0$$ which together
with an application of Alexander duality gives $H_1(M_i^3)\cong H_1(M_i^3,
M_i^3 - \Sigma_i)\cong H^2(\Sigma_i)$.  Further, as the isomorphisms
$H_1(M_i^3)\cong H^2(M_i^3)$ and $H_1(M_i^3, M_i^3 - \Sigma_i)\cong
H^2(\Sigma_i)$, given respectively by Poincar\'e and Alexander duality,
are obtained by taking cup products with the fundamental class, the
diagram
\[\begin{array}{ccc}
H^2(M_i^3)   & \longrightarrow & H^2(\Sigma_i)\\
 \downarrow &  &  \downarrow \\
 H_1(M_i^3)  &\longrightarrow & H_1(M_i^3, M_i^3 - \Sigma_i)
\end{array}\]
commutes. 

\vspace{0.2cm}
\noindent 
Thus the inclusion of $\Sigma_i$ in $M_i^3$ gives an isomorphism
$H^2(M_i^3)\cong H^2(\Sigma_i)$. Since $H_2(M_i^3)$ and $H_2(\Sigma_i)$
have no torsion, the cap product induces perfect pairings
$H_2(M_i^3)\times H^2(M_i^3)\to \mathbb{Z}$ and $H_2(\Sigma_i)\times
H^2(\Sigma_i)\to \mathbb{Z}$. Therefore, by duality, 
the map $H_2(\Sigma_i)\to H_2(M_i^3)$
induced by inclusion is also an isomorphism.

\vspace{0.2cm}
\noindent 
Now take a basis for $H^2(\Sigma_i)$. Each element of this basis can be
looked at as an integral linear combination of the planar surfaces (as
in cellular homology), with trivial boundary. We obtain a surface
corresponding to each such homology class by taking copies of the
planar surfaces, with the number and orientation determined by the
coefficient. Since the homology classes are cycles, these planar
surfaces can be glued together at the boundaries to form closed,
oriented, immersed surfaces. Without loss of generality, we can assume
these to be connected.

\begin{remark} By doing surgeries on the seams of $\Sigma_i\subset M_i^3$
some new homology is created (the homology of $M_{i+1}^3$)
One constructs naturally surfaces representing the homology  of
$M_{i+1}^3$, as follows. 
One considers {\em generalized} Seifert surfaces in $M_i^3$ of
the loci of the surgeries, which are surfaces which might 
have boundary components along other seams. Then one caps-off 
the boundaries by using the cores of the 2-handles which are added and
push the closed surfaces into $M_{i+1}^3$. 
Notice that we can consider also some Seifert surface whose boundary 
components are seams in some $M_{i+n}^3$ for $n > 1$. 
\end{remark}

\subsection{On Casson finiteness}\label{S:proof}

\noindent 
Suppose we do have a $4$-manifold $(W^4,\del W^4)$ with a g.s.c.
handle decomposition of its interior. Since there
may be infinitely many handles, we cannot use Casson's
argument. However, we note that we can use Casson's argument if we
can show that
\begin{itemize}
\item $(W^4,\del W^4)$ has a (finite) handle decomposition without
$1$-handles.
\item Some $(Z^4,\del Z^4)$ has a handle decomposition without
$1$-handles, where $Z^4$ is compact, contractible with $\pi_1(\del
Z^4)=\pi_1(\del W^4)$.
\item Some $(Z^4,\del Z^4)$ has a handle decomposition without
$1$-handles, where $Z^4$ is compact, contractible and there is a
surjection $\pi_1(\del Z^4)\onto \pi_1(\del W^4)$ (by
Proposition~\ref{T:alg}).
\end{itemize}
\noindent Thus, we can apply Casson's argument if we show finiteness, or some
weak form of finiteness such as the latter statements above.

\vspace{0.2cm}
\noindent 
We now assume that the handle decomposition is as in the conclusion of
Theorem~\ref{T:reorder}.  We will change our measures of time so that
passing from $M_i^3$ to $M_{i+1}^3$ consists of performing all the
surgeries required to compress the surfaces, splitting along the
2-spheres, and also performing the necessary surgeries off the surface.

\vspace{0.2cm}
\noindent 
In $M_i^3$, we have a collection of embedded surfaces representing all
the homology of $M_i^3$. We see that we have Casson finiteness in a
special case.

\begin{theorem}\label{embed} If there exists  $i$ so that  the immersion 
$\psi_i:\cup_j F_j^2(i)\to M_i^3$ is actually an embedding and 
$\phi_i(\psi_i(F_j^2(i))=\{1\}\subset\pi_1(\del W^4)$ 
then $\pi_1(\del W^4)$ violates the Kervaire conjecture.
\end{theorem}
\begin{proof} Let $k$ be the rank of $H_1(M_i^3)$ and $P_j,1\le j\le k$ be
the fundamental groups of the surfaces. Since the surfaces are
disjoint, $\pi_1(M_i^3)$ is obtained by HNN extensions from the
fundamental group $G$ of the complement of the surfaces. Thus, if
$\psi_j$ are the gluing maps, we have
        $$\pi_1(M_i^3)=\left<G,t_1,\dots,t_k;t_j x
t_j^{-1}=\psi_j(x)\forall x\in P_j\right>$$ 
Now, since $\phi_i(P_j)=1$ and $\phi_i(G)=\pi_1(\del
W^4)$, $\pi_1(M_i^3)$ surjects onto $\left<\pi_1(\del
W^4),t_1,\dots,t_n\right>$, the group obtained by adding $k$ generators
to $\pi_1(\del W^4)$. But, $M_i^3$ is obtained by using $n$ 2-handles and
$n-k$ 3-handles. Thus, as in Casson's theorem, $\pi_1(M_i^3)$ is killed
by adding $n-k$ generators and $n$ relations. 
This implies that $\pi_1(\del W^4)$ is killed
by adding $n$ generators and $n$ relations.
\end{proof}

\begin{theorem}\label{imers} There exists always an immersion 
as in the structure theorem with the additional property that 
$\phi_i(\psi_i(F_j^2(i))=\{1\}\subset\pi_1(\del W^4)$ holds true. 
\end{theorem}
\begin{proof}
By construction the images of the seams (which are roughly 
speaking half the generators of the fundamental group) are null-homotopic. 
If the fundamental groups of the 
generalized  Seifert surfaces from the previous remark 
map to the trivial group, then  after 
doing surgery on the seams  we
obtain surfaces $F_j^2(i+1)$ representing homology with trivial $\pi_1$
images by $\phi_{i+1}$. 
Thus, it suffices to show that we obtain this condition for a choice 
of Seifert surfaces for all seams, at some time in the future.

\vspace{0.2cm}
\noindent 
Fix $i$ large enough so that $M_i^3$ is in the collar 
$\del W^4 \times[0,\infty)$. 
\begin{lemma}\label{n'} There exists some $n'=n'(i)$ such that, whenever 
a 2-sphere immersed in ${\rm int}(W^4)-K_{i+n'}^4$ bound a 3-ball 
immersed in the collar, it actually bounds a 3-ball immersed 
in ${\rm int}(W^4)-K_i^4$. 
\end{lemma}
\begin{proof}
Choose $n'$ large enough so that a small collar 
$\del W^4\times [x,y]\subset K_{i+n'}^4-{\rm int}(K_{i}^4)$.
Then use the horizontal flow to send 
$\del W^4\times [0,y]$ into $\del W^4\times [x,y]$ by preserving the 
right side boundary. This yields a ball in the complement of $K_i^4$. 
\end{proof}

\noindent We need the following analogue, for 
$\pi_2$ instead of $\pi_1$,  of Lemma~\ref{T:persistsdies}. 
\begin{lemma}
There exists $k=k(i)$ such that any immersed 2-sphere in $K_i^4$ 
(respectively in the intersection of $K_i^4$ with the collar) 
which bounds a 3-ball in ${\rm int}(W^4)$ (respectively in the collar) 
 does so in $K_{i+k}^4$ (respectively  in the collar).
\end{lemma}  

\begin{proof}
The same trick we used in the proof of the previous lemma applies. 
\end{proof}

\noindent 
Choose now $n$ large enough such that 
$K_{i+n}^4-K_{i+n'}^4$ contains a non-trivial collar $\del W^4\times [x,y]$,
and $n > k(i+n')$ provided by the Lemma~\ref{T:persistsdies}. 
Consider a surgery curve $\gamma\subset M_j^3$ for some 
$i+n' \leq j < i+n$. 

\begin{lemma}
There exists a generalized Seifert surface for $\gamma$ in $M_i^3$, 
so that the other boundary components are surgery curves 
from $M_m^3$, with $m\leq i+n$, and 
whose fundamental group maps to the trivial group under $\phi_i$. 
\end{lemma}

\begin{proof}
The curve $\gamma$ bounds a disc in the $2$-handle attached to
it. Further, as it can be pulled back, say along an annulus, to time
$M_{i+n'}^3$, and then dies by $M_{i+n}^3$, it bounds another disc consisting
of the annulus and the disc by which it dies. These discs together
form an immersed  $2$-sphere. Consider the class 
$\nu\in \pi_2(\del W^4)$ of this 2-sphere
by using the projection of the collar 
on $\del W^4$. We can  realize the element $\nu$ by an immersed 
2-sphere in a small collar $\del W^4\times [x,y]\subset K_{i+n}^4-K_i^4$. 
Therefore by modifying the initial 2-sphere by this sphere
(which is far from $\gamma$) in the small collar 
one finds an immersed 2-sphere whose image in  $\pi_2(\del W^4)$
is trivial. Since $i$ was large enough $K_{i+n}^4-K_i^4$ is a subset 
of a larger collar $\del W^4 \times [0, z]$. Then the 2-sphere we 
constructed bounds a 3-ball in $\del W^4 \times [0, z]$, and so by 
Lemma~\ref{n'} it also does so in $X^4=K_{i+n}^4-K_i^4$. 
Let $\mu:S^2\to X^4$ denote this immersion
realizing a trivial element of $\pi_2(X)$. Then $\mu$ lifts to a map  
$\tilde \mu: S^2\to \tilde X^4$, where $\tilde X^4$ is the universal 
covering space of $X^4$. 
Since $\mu$ is null-homotopic 
the homology  class of $[\tilde\mu]=0\in H_2(\tilde X^4)$ 
is trivial, when interpreting $\tilde \mu$ as a 2-cycle in $\tilde X^4$. 

\vspace{0.2cm}
\noindent 
The homology of $\tilde X^4$ is computed from the 
$\pi_1(X^4)$-equivariant complex associated to the handle
decomposition, whose generators in degree $d$ are the $d$-handles 
attached to $K_i^4$ in order to get $K_{i+n}^4$. 
Therefore one has then the following relation in this differential complex:  
\[ [\tilde\mu]= \sum_j c_j\, \partial\, [h_j^3], \;\; c_j\in \mathbb Z. \]
The action of the algebraic boundary operator $\partial$ on 
the element  $[h_j^3]$ can be described  
geometrically as the class of the 2-cycle which represents the
attachment 2-sphere $\partial^+ h_j^3$ of the 3-handle $h_j^3$. 
Consequently the previous formula can be rewritten as 
\[ \tilde\mu = \sum_j c_j \,\partial^+ h_j^3 + \sum _k d_k \,L_k, \;\;
c_j, d_k \in \mathbb Z\]
where $L_k$ are closed  surfaces (actually these are closed 2-cycles,
but they can be represented by surfaces by the well-known results of
R.~Thom)  with the property that 
\[ [L_k]\cdot [\delta_m^2]=0, \forall k,m\]
($\delta_m^2$ denotes the core of the 2-handle $h_m^2$).
Let us compute explicitly the boundary operator on the 3-handles, in
terms of the surfaces we have in the 2-complex $\Sigma_i$. 
Set
\[ \partial\, [h_j^3]= \sum_k m_{jk}\, [h_k^2].\]
Then the  coefficient $m_{jk}$ is the number of times 
the boundary $\partial^+h_j^3$ runs over the core of $h_k^2$. 
But the 2-sphere $\partial^+h_j^3$, when pulled back in $M_i^3$, is a 
planar surface in $M_i^3$ whose boundary circles (i.e. at seams) 
are capped-off by the core disks $\delta_k^2$ of the 2-handles $h_k^2$. 
Therefore the number $m_{jk}$ is the number of times the seam 
$\partial^+h_k^2$ appears in the planar surface which is a pull-back of
$\partial^+h_j^3$.
In particular the coefficient of a $2$-handle vanishes in a
$3$-cycle only if the boundaries of the planar surface glue 
together to close up at the corresponding surgery locus.

\vspace{0.2cm}
\noindent 
Thus the pull-backs of the surfaces $\sum_j c_j \partial^+ h_j^3$
give a surface $F^2$ in $M_i^3$ with boundary the
curve $\gamma$ with which we started, plus some other 
curves along which surgery is performed by time $i+n$. 
As this is in fact a closed cycle in the
universal cover, the surface $F^2$ lifts to a surface in  $\tilde X^4$, 
with a single boundary component, corresponding to 
a curve which is not surgered by the time $i+n$.
Therefore the map $\pi_1(F^2)\to \pi_1(X^4)\to \pi_1(\partial W^4)$ 
factors through $\pi_1(\tilde X^4)=1$, hence  the
image of $\pi_1(F^2)$ in $\pi_1(\del W^4)$ is trivial.
\end{proof}
\noindent 
Thus, after surgering along the curves up to the $M_{i+n}^3$  
we do have the required Seifert surfaces to compress to get
embedded surfaces with trivial $\pi_1(\del W^4)$ image.
\end{proof}

\begin{proposition}
If $\gamma_{i,k}\subset M_i^3$ are the surgery curves (i.e. the seams) 
then $\gamma_{i,k}\in LCS_{\infty}(\pi_1(M_i^3))$, where 
$LCS_s(G)$ is the lower central series of the group
$G$, $LCS_1(G)=G, LCS_{s+1}(G)=[G,LCS_{s}(G)]$, and 
$LCS_{\infty}(G)=\cap_{s=1}^{\infty}LCS_s(G)$. 
\end{proposition}
\begin{proof} We will express each surgery locus $\gamma$ as a product of
commutators of the form $[\alpha_i,\beta_i]$, with each $\alpha_i$
being conjugate to a surgery locus (possibly $\gamma$ itself). It then
follows readily that $\gamma\in LCS_{\infty}(\pi_1(M_i^3))$, as now if each
$\gamma_{i,k}\in LCS_s(\pi_1(M_i^3))$, 
then each $\gamma_{i,k}\in LCS_{s+1}(\pi_1(M_i^3))$.

\vspace{0.2cm}
\noindent 
Suppose now $\gamma=\gamma_{i,k}$ is a surgery locus. Then the
$0$-frame surgery along $\gamma$ creates homology in $M_{i+1}^3$ which
by our structure theorem is represented by a surface
$S^2=\psi_{i+1}(F_j^2(i+1))$. The pullback of $S^2$ to $i$ gives a surface
with boundary along seams, and being compressed to a sphere by the
seams, so that the algebraic multiplicity of $\gamma$ is $1$ while
that of all other seams is $0$. In terms of the fundamental group,
this translates to the relation that was claimed.
\end{proof}

\noindent 
Since, we have immersed surfaces of the required form, the obstruction
we encounter is in making these surfaces disjoint at some finite
stage. Note that for a finite decomposition, we do indeed have
disjoint surfaces representing the homology after finitely many
surgeries, since we in fact have a family of such spheres.

\begin{figure}
\hspace{3cm}\includegraphics{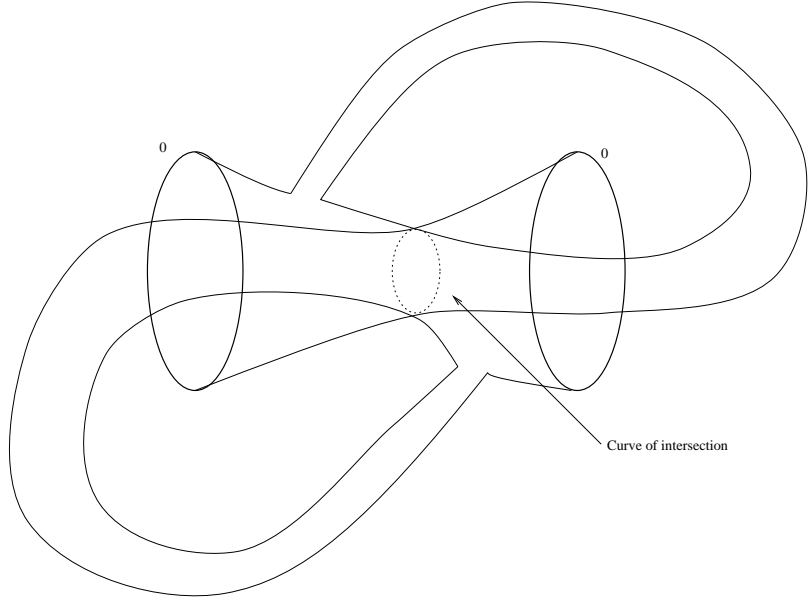}
\caption{Tori with one intersection curve}\label{F:plumb}
\end{figure}

\begin{proposition} There exists a $2$-complex 
$\Sigma=\cup_{i=1}^{\infty}\Sigma_i$ with intersections
along double curves, coming from a handle-decomposition as above,
where all the seams are trivial in homology, but which does not carry
disjoint, embedded surfaces representing all of the homology.
\end{proposition}
\begin{proof} For the first stage, take two surfaces of genus $2$, and
let them intersect transversely along two curves (which we  call
seams) that are disjoint and homologically independent in each
surface. Next, take as Seifert surfaces for these curves once
punctured surfaces of genus $2$ intersecting in a similar manner, and
glue their boundary to the above-mentioned curves of
intersection. Repeat this process to obtain the complex.

\vspace{0.2cm}
\noindent 
At the first stage, we cannot have embedded, disjoint surfaces
representing the homology as the cup product of the surfaces is
non-trivial. As the surfaces are compact, we must terminate at some
finite stage. We will prove that if we can have disjoint surfaces at
the stage $k+1$, then we do at stage $k$. This will suffice to give
the contradiction.

\vspace{0.2cm}
\noindent 
Now, we know the complex cannot be embedded in the first
stage. Suppose we did have disjoint embedded surface $F_1$ and $F_2$
at stage $k+1$. Since these form a basis for the homology, they
contain curves on them that are the seams at the first stage with
algebraically non-zero multiplicity, i.e., the collection of curves
representing the seam is not homologically trivial in the intersection
of the first stage with the surface. Further, some copy of the first
seam must bound a subsurface $F_i'$ in each of the surfaces, for
otherwise the surface contains a curve dual to the seam. For, the cup
product of such a dual curve with the homology class of the other
surface is non-trivial, hence it must intersect the other surface,
contradicting the hypothesis that the surfaces are
disjoint. Similarly, at the other seam we get surfaces $F_i''$.

\vspace{0.2cm}
\noindent 
By deleting the first stage surfaces and capping off the first stage
seams by attaching discs, we get a complex exactly as before with the
$(j+1)$th stage having become the $j$th stage. Further, the $F_1'$ and
$F_2''$ now give disjoint, embedded surfaces representing the homology
that are supported by stages up to $k$. This suffices as above to
complete the induction argument.

\vspace{0.2cm}
\noindent 
It is easy to construct a handle-decomposition corresponding to this
complex. Figure~\ref{F:plumb} shows a construction of tori with one
curve of intersection. Here we have used the notation of Kirby
calculus, with the thickened curves being an unlink along each
component of which $0$-frame surgery has been performed. It is easy to
see that the same construction can give surfaces of genus $2$
intersecting in $2$ curves. On attaching the first two $2$-handles,
the boundary is $(S^2\times S^1)\# (S^2\times S^1)$. Since the curves
of intersection are unknots, after surgery they bound
spheres. Further, it is easy to see by cutting along these that the
boundary is $(S^2\times S^1)\# (S^2\times S^1)$ after attaching the
$2$-handles and $3$-handles as well. Repeating this process, we obtain
our embedding.

\vspace{0.2cm}
\noindent 
Thus we have an infinite handle-decomposition satisfying our
hypothesis for which this $2$-complex is $\Sigma$.
\end{proof}

\subsection{A wild example}

\noindent 
We will construct an example of an open, contractible $4$-manifold
that is not tame, and that has a handle-decomposition without
$1$-handles.

\begin{theorem} There is a proper handle-decomposition of an open,
contractible $4$-manifold $W^4$ such that $W^4$ is not the interior of a
compact $4$-manifold. In particular $W^4$ does not have a finite
handle-decomposition.
\end{theorem}
\begin{proof} We will take a variant of the example in the last
section. Namely, we construct an explicit handle-decomposition
according to a canonical form.

\vspace{0.2cm}
\noindent 
Start with a $0$-handle and attach to its boundary three $2$-handles
along an unlink. the resulting manifold has boundary $(S^2\times S^1)\#
(S^2\times S^1)\#(S^2\times S^1)$ obtained by $0$-frame-surgery about
each component of an unlink with $3$ components. We now take as
Seifert surfaces for these components surfaces of genus $2$, so that
each pair intersects in a single curve, so that the curves of
intersection form an unlink and are unlinked with the original
curves. 

\vspace{0.2cm}
\noindent 
Now, attach $2$-handles along the curves of intersections, and then
$3$-handles along the Seifert surfaces compressed to spheres by adding
discs in the $2$-handles just attached. It is easy to see that the
resulting manifold once more has boundary $(S^2\times S^1)\# (S^2\times
S^1)\#(S^2\times S^1)$. Thus, we may iterate this process. Further, the
generators of the fundamental group at any stage are the commutators
of the generators at the previous stage.

\vspace{0.2cm}
\noindent 
Suppose $W^4$ is in fact tame. Then, we may use the results of the
previous sections. Now, by construction no curve dies as only trivial
relations have been added. Thus every element in kernel($\phi_i$) must
fail to persist by some uniform time. In particular, the image of the
group after that time in the present (curves that persist beyond that
time) must inject under $\phi_i$. But we know that it also
surjects. Thus, we must have an isomorphism.

\vspace{0.2cm}
\noindent 
Thus, there is a unique element mapping onto each element of
$\pi_1(\del W^4)$. Hence this element must persist till infinity as we
have a surjection at all times. On the other hand, since the limit of
the lower central series of the free group is trivial, no non-trivial
element persists. This gives a contradiction unless $\pi_1(\del W^4)$ is
trivial.

\vspace{0.2cm}
\noindent 
But there are non-trivial elements that do persist beyond any give
time. As no element dies, we again get a contradiction.
\end{proof}

\subsection{Further obstructions from Gauge theory}\label{S:gauge}

\noindent 
To further explore some of the subtleties that one might encounter in
trying to construct a handle decomposition without $1$-handles for a
contractible manifold, given one for its interior, we consider a more
general situation. We will consider sequences of $3$-manifolds $M_i^3$
that begin with $S^3$. As before, we require that each manifold comes
from the previous one by $0$-frame surgery about a homologically
trivial curve or by splitting along a non-separating $S^2$ and capping
off. Also, we require degree-one maps $f_i$ to a common manifold $N^3$,
related as before. We will say that the sequence limits to $N^3$ if
{\em any curve that persists dies} as in lemma~\ref{T:persistsdies}.

\vspace{0.2cm}
\noindent 
In this situation, our main question generalizes to a {\em relative
version}, namely, given any such sequence $\{M_i^3\}$, with $M_k^3$ an
element in the sequence, is there a finite sequence that agrees up to
$M_k^3$ with the old sequence and whose final term is $N^3$?

\vspace{0.2cm}
\noindent 
We will show that there is an obstruction to completing certain
sequences to finite sequences when $N^3=S^3$. We do not know whether
there are infinite sequences limiting to $N^3$ in this case.

\vspace{0.2cm}
\noindent 
Let $\P^3$ denote the Poincar\'e homology sphere. Observe that we cannot
pass from this to $S^3$ by $0$-frame surgery about homologically
trivial curves and capping-off non-separating spheres. For, if we
could, $\P^3$ would bound a manifold with $H^2=\oplus_k\H$, which is
impossible as $\P^3$ has Rochlin invariant $1$. On the other hand, for
the same reason, $\P^3$ cannot be part of any sequence of the above
form.

\vspace{0.2cm}
\noindent 
Using Donaldson's theorem~\cite {Do}, we have a similar result for the
connected sum $\P^3\#\P^3$ of $\P^3$ with itself. The main part of the proof
of this lemma was communicated to us by R.~Gompf.

\begin{lemma}\label{T:PP} One cannot pass from $\P^3\#\P^3$ to $S^3$ by
$0$-frame surgery along homologically trivial curves and capping off
non-separating $S^2$'s.
\end{lemma}
\begin{proof} If we did have such a sequence of surgeries, then
$\P^3\#\P^3$ bounds a $4$-manifold $M_1^4$ with $H^2=\oplus_k\H$ , with a
half-basis formed by embedded spheres. Now glue this to a manifold
with form $E_8\oplus E_8$ which is bounded by $\P^3\#\P^3$ to get $M^4$.

\vspace{0.2cm}
\noindent 
We can surger out the disjoint family of $S^2$'s from $M^4$ to get a
$4$-manifold with form $E_8\oplus E_8$ and trivial $H_1$. This
contradicts Donaldson's theorem.
\end{proof}

\noindent 
We still do not have a sequence as claimed. For, Cassson's argument
shows that $\P^3\#\P^3$ cannot be part of a sequence. To obtain such a
sequence, we will construct a manifold $N^3$ that can be obtained by
$0$-frame surgery on algebraically unlinked $2$-handles from each of
$S^3$ and $\P^3\#\P^3$. Thus, $N^3$ is part of a sequence. On the other
hand, if we had a sequence starting at $N^3$ that terminated at $S^3$,
then we would have one starting at $\P^3\#\P^3$, which contradicts the
above lemma.

\vspace{0.2cm}
\noindent 
To construct $N^3$, take a contractible $4$-manifold $K^4$ that bounds
$\P^3\#\P^3$. By Freedman's theorem (\cite {Fr}), this exists, and can
moreover be smoothed after taking connected sums with sufficiently
many copies of $S^2\times S^2$. Take a handle decomposition of
$K^4$. This may include $1$-handles, but these must be boundaries of
$2$-handles. Hence, by handle-slides, we can ensure that each
$1$-handle is, at the homological level, a boundary of a $2$-handle
and is not part of the boundary of any other $2$-handle. Replacing the
$1$-handle by a $2$-handle does not change the boundary, and changes
$H^2(K^4)$ to $H^2(K^4)\oplus \left(\oplus_k\H \right)$. We do this dually
with $3$-handles too. Sliding $2$-handles over the new ones, we can
ensure that the attaching maps of the $2$-handles are having the same
algebraic linking (and framing) structure as a disjoint union of Hopf
links.

\vspace{0.2cm}
\noindent 
Now take $N^3$ obtained from $S^3$ by attaching half the links, so that
these are pairwise algebraically unlinked. The manifold $N^3$ has the
required properties.


\begin{thebibliography}{}

\bibitem{Ba}
G.~Baumslag,
\newblock Topics in Combinatorial Group Theory, 
\newblock {\em Lectures in Math., ETH Zurich, Birkhauser}, 1993.

\bibitem{BM}
S.G.~Brick and M.L.~Mihalik,
\newblock {\em The QSF property for groups and spaces},
\newblock  {Math.Zeitschrift} {\bf 220}(1995), 207--217.

\bibitem{BT}
M.G.~Brin and T.L.~Thickstun, 
{\em $3$-manifolds which are end $1$-movable}, Mem. A.M.S. {\bf 81}(1989),
no. 411, 73p.


\bibitem{BZ}
M.~Boileau and H.~Zieschang,
\newblock {\em Heegaard genus of closed orientable 
Seifert 3-manifolds},
\newblock {Inventiones Math.} {\bf 76}(1984), 455-468.


\bibitem{Do}
S.K.~Donaldson,
\newblock \emph{An application of gauge theory to $4$-dimensional topology},
\newblock J.Diff.Geometry  {\bf 26}(1983), 279--315.

\bibitem{Fr}
M.~Freedman, 
\newblock {\em The topology of four-dimensional manifolds},
\newblock {J.Diff.Geometry} {\bf 17}(1982), 357-453.


\bibitem{Fun2}
L.~Funar,
\newblock {\em On proper  homotopy type and the simple connectivity at 
infinity of open 3-manifolds}, 
\newblock Atti Sem.Mat.Fis. Univ. Modena {\bf 49}(2001), 15-29.



\bibitem{FT}
L.~Funar and T.L.~Thickstun,
\newblock {\em On open 3-manifolds proper homotopy equivalent to 
geometrically simply connected polyhedra}, 
\newblock {Topology Appl.} {\bf 109}(2001), 191-200.

\bibitem{GR}
M.~Gerstenhaber and O.S.~Rothaus,
\emph{The solution of sets of equations in groups},
Proc. Nat.Acad.Sci. USA  {\bf 48}(1962), 1531--1533.


\bibitem{Glaser}
L.C.~Glaser,
\newblock {\em Uncountably many contractible 4-manifolds},
\newblock {Topology} {\bf 6}(1965), 37-42.


\bibitem{HP}
A.~Haefliger and V.~Po\'enaru,
\newblock {\em La clasification des immersions combinatoires},
\newblock Publ.Math. I.H.E.S. {\bf 23}(1964), 75-91. 


\bibitem{HR}
B.~Hughes and A.~Ranicki, 
\newblock {Ends of complexes}, 
\newblock {\em Cambridge Univ.Press}, 123, 1996.  


\bibitem{Ke}
M.~Kervaire, 
\newblock \emph{On higher dimensional knots},
\newblock Differential and combinatorial topology (Symposium in honour
of Marston Morse), Princeton Math. Series, vol. 27, 1965.

\bibitem{LiS}
W.B.R.~Lickorish and L.C.~Siebenmann, 
\newblock {\em Regular neighborhoods and the stable range},
\newblock {Trans. A.M.S.} {\bf 139}(1969), 207-230. 


\bibitem{LS}
R.C.~Lyndon and P.E.~Schupp,
\newblock {Combinatorial Group Theory},
\newblock {\em Ergebnisse Math., Springer-Verlag}, 89, 1977.


\bibitem{MKS}
W.~Magnus, A.~Karras and D.~Solitar,
{Combinatorial Group Theory},
\emph{ Interscience}, 1966.

\bibitem{Mazur0}
B.~Mazur,
\newblock {\em A note on some contractible 4-manifolds},
\newblock {Ann. of Math.} {\bf 73}(1961), 221-228. 

\bibitem{Mazur}
B.~Mazur,
\newblock {\em Differential topology from the point of view of simple
  homotopy theory}, 
\newblock {Publ.Math. I.H.E.S.} {\bf 15}(1963), 1--93.

\bibitem{McM}
D.R.~McMillan Jr,
\newblock {\em Some contractible  open 3-manifolds},
\newblock {Trans. A.M.S.} {\bf 102}(1962), 373-382. 

\bibitem{M}
M.~Mihalik, 
{\em Compactifying coverings of $3$-manifolds}, 
Comment.Math.Helv. {\bf 71}(1996), 362--372.

\bibitem{Po}
V.~Po\'enaru,
\newblock {\em La d\'ecomposition de l'hypercube en produit topologique},
\newblock {Bull. Soc. Math. France} {\bf 88}(1960), 113-129. 


\bibitem{Po0}
V.~Po\'enaru.
\newblock {\em On the equivalence relation forced by the singularities of a non-degenerate simplicial map},
\newblock {Duke Math.J.}  {\bf 63}(1991), 421--429.


\bibitem{Po1}
V.~Po\'enaru,
\newblock {\em Killing handles of index one stably and
  $\pi_1^{\infty}$}, 
\newblock {Duke Math.J.}  {\bf 63}(1991), 431--447.


\bibitem{Po4}
V.~Po\'enaru,
\newblock {\em $\pi_1^{\infty}$ and infinite simple homotopy type in dimension three},
\newblock {in "Low Dimensional Topology", Contemporary
     Math., {\bf 233}(1999), (H.Nencka, Editor), AMS, 
     1-28}. 

\bibitem{PoU}
V.~Po\'enaru,
\newblock {Universal covering spaces of closed 3-manifolds are simply-connected
at infinity}, 
\newblock preprint Orsay, 20/2000.


\bibitem{PoTa}
V.~Po\'enaru and  C.Tanasi,
\newblock {\em Some remarks on geometric simple connectivity},
\newblock {Acta Math. Hungar.} {\bf 81}(1998), 1-12.

\bibitem{Quinn}
F.~Quinn,
\newblock {\em The stable topology of 4-manifolds},
\newblock {Topology Appl.}  {\bf 15}(1983), 71-77.


\bibitem{Rot}
O.S.~Rotthaus, 
\newblock {\em On the non-triviality of some group extensions 
given by generators and relations},
\newblock Ann. of Math. {\bf 106}(1977), 599-612.

\bibitem{RS}
C.P.~Rourke and B.J.~Sanderson,
Introduction to PL topology, 
{\em Springer Verlag}, 1972.

\bibitem{Si}
L.~Siebenmann, 
\newblock {\em The obstruction to finding a boundary for an open 
manifold of dimension greater than five}, 
\newblock {PhD Thesis, Princeton}, 1965. 

\bibitem{Th}
W.P.~Thurston,
\newblock {\em Three dimensional manifolds, Kleinian groups and 
hyperbolic geometry},
\newblock {Bull. A.M.S.} {\bf 6}(1982), 357-381.


\bibitem{Wall}
C.T.C.~Wall, 
\newblock {\em Geometrical connectivity I, II},
\newblock J.London Math.Soc. (2) {\bf 3}(1971), 597-604, 605-608.

\bibitem{Wr}
D.G.~Wright,
\newblock {\em Contractible open manifolds which are not covering spaces},
\newblock {Topology} {\bf 31}(1992), 281-291.

\bibitem{Ze}
E.C.~Zeeman,
\newblock {\em The Poincar\'e conjecture for $n\geq 5$},
\newblock Topology of 3-manifolds and related topics, 
(M.K.Fort Jr., Editor), Prentice Hall, N.J., 1962, 198-204.

\end{thebibliography}
\end{document}